\documentclass[12pt]{amsart}
\usepackage{tikz}
\usetikzlibrary{matrix, calc, arrows,backgrounds}

\usepackage{float}
\tikzset{node distance=2cm,auto}
\usepackage{amsfonts,amsmath,amssymb,amsthm,mathrsfs,verbatim}
\usepackage[all,cmtip]{xy}
\usepackage{latexsym}
\usepackage{graphics}
\newtheorem{theorem}{Theorem}[section]
\newtheorem{lemma}[theorem]{Lemma}
\newtheorem{corollary}[theorem]{Corollary}
\newtheorem{proposition}[theorem]{Proposition}
\newtheorem{conjecture}{Conjecture}

\newtheorem{mainthm}{Theorem}

\newtheorem{maincor}[mainthm]{Corollary}

\theoremstyle{definition}
\newtheorem{question}[theorem]{Question}
\newtheorem{definition}[theorem]{Definition}

\newtheorem{assumption}{Assumption}
\newtheorem{example}[theorem]{Example}

\newtheorem{remark}[theorem]{Remark}

\usepackage{amsfonts}
\usepackage{color}

\newcommand{\bpf}{\noindent{\bf Proof}\hspace{7pt}}
\newcommand{\epf}{\qed}
\newcommand{\ben}{\begin{enumerate}}
\newcommand{\een}{\end{enumerate}}
\newcommand{\ble}{\begin{lemma}}
\newcommand{\ele}{\end{lemma}}
\newcommand{\bth}{\begin{theorem}}
\renewcommand{\eth}{\end{theorem}}
\newcommand{\bmth}{\begin{mainthm}}
\newcommand{\emth}{\end{mainthm}}
\newcommand{\bmco}{\begin{maincor}}
\newcommand{\emco}{\end{maincor}}

\newcommand{\bpr}{\begin{proposition}}
\newcommand{\epr}{\end{proposition}}
\newcommand{\bco}{\begin{corollary}}
\newcommand{\eco}{\end{corollary}}
\newcommand{\bcon}{\begin{conjecture}}
\newcommand{\econ}{\end{conjecture}}
\newcommand{\bqu}{\begin{question}}
\newcommand{\equ}{\end{question}}
\newcommand{\bde}{\begin{definition}}
\newcommand{\ede}{\end{definition}}
\newcommand{\bas}{\begin{assumption}}
\newcommand{\eas}{\end{assumption}}
\newcommand{\bre}{\begin{remark}}
\newcommand{\ere}{\end{remark}}
\newcommand{\bex}{\begin{example}}
\newcommand{\eex}{\end{example}}
\newcommand{\barr}{\begin{array}}
\newcommand{\earr}{\end{array}}
\newcommand{\btab}{\begin{tabular}}
\newcommand{\etab}{\end{tabular}}
\newcommand{\beq}{\begin{equation}}
\newcommand{\eeq}{\end{equation}}
\newcommand{\bea}{\begin{eqnarray*}}
\newcommand{\eea}{\end{eqnarray*}}
\newcommand{\bce}{\begin{center}}
\newcommand{\ece}{\end{center}}
\newcommand{\bpi}{\begin{picture}}
\newcommand{\epi}{\end{picture}}
\newcommand{\bfi}{\begin{figure} \begin{center}}
\newcommand{\efi}{\end{center} \end{figure}}

\newcommand{\bsl}{\begin{slide}{}}
\newcommand{\esl}{\end{slide}}

\newcommand{\ul}{\underline}

\newcommand{\hso}[1]{\hspace{-1pt}}

\newcommand{\into}{\hookrightarrow}

\newcommand{\sbe}{\subseteq}

\def\<{\langle}
\def\>{\rangle}

\newcommand{\cD}{{\mathcal D}}

\newcommand{\cE}{{\mathcal E}}

\newcommand{\cH}{{\mathcal H}}

\newcommand{\cP}{{\mathcal P}}

\newcommand{\cQ}{{\mathcal Q}}

\newcommand{\cS}{{\mathcal S}}

\newcommand{\cU}{{\mathcal U}}

\newcommand{\sCT}{{\mathsf{CT}}}
\newcommand{\sPh}{{\mathsf{Ph}}}

\newcommand{\twA}{{}^2\! {A}}
\newcommand{\twD}{{}^2\! {D}}

\newcommand{\twE}{{}^2\! {E}}

\newcommand{\sfE}{\mathsf E}

\newcommand{\Aut}{\mathop{\rm Aut}\nolimits}

\renewcommand{\bar}{\overline}

\newcommand{\diag}{\mathop{\rm diag}\nolimits}

\DeclareMathOperator{\id}{id}
\newcommand{\im}{\mathop{\rm im}\nolimits}

\DeclareMathOperator{\Syl}{Syl}

\DeclareMathOperator{\SL}{SL}
\DeclareMathOperator{\PSL}{PSL}

\def\cD{{\mathcal D}}

\newcommand{\FF}{{\mathbb F}}

\newcommand{\NN}{{\mathbb N}}

\newcommand{\ZZ}{{\mathbb Z}}

\newcommand{\hform}{{\sf h}}

\newcommand{\dfn}{\em}

\newcommand{\after}{\mathbin{ \circ }}
\DeclareMathOperator{\GL}{GL}
\DeclareMathOperator{\PGL}{PGL}
\DeclareMathOperator{\GD}{T}

\DeclareMathOperator{\SU}{SU}

\DeclareMathOperator{\Sp}{Sp}
\DeclareMathOperator{\PSp}{PSp}
\DeclareMathOperator{\GU}{GU}
\DeclareMathOperator{\PGU}{PGU}
\DeclareMathOperator{\GSp}{GSp}

\newcommand{\GamL}{\mathop{\rm \Gamma L}}

\DeclareMathOperator{\PGammaL}{P\Gamma L}
\DeclareMathOperator{\GammaSp}{\Gamma Sp}
\DeclareMathOperator{\PGammaSp}{P\Gamma Sp}
\DeclareMathOperator{\GammaU}{\Gamma U}
\DeclareMathOperator{\PGammaU}{P\Gamma U}
\DeclareMathOperator{\PGSp}{PGSp}

\newcommand{\Stab}{\mathop{\rm Stab}}

\newcommand{\normal}{\lhd}
\newcommand{\mspl}{.}

\newcommand{\Chi}{{\mathcal X}}

\renewcommand{\hat}{\widehat}

\newcommand{\vep}{\varepsilon}

\newcommand{\trin}{\tau}

\renewcommand{\qed}{\hfill $\square$}

\newcounter{romanlistctr}
{\end{list}}%

 \makeatletter
 \@addtoreset{equation}{section}
 \makeatother

 \makeatletter
 \def\section{\@startsection {section}{1}{\z@}{-1.5ex plus -.5ex
 minus -.2ex}{1ex plus .2ex}{\large\bf}}

 \makeatletter
 \def\subsection{\@startsection {subsection}{1}{\z@}{-1.5ex plus -.5ex
 minus -.2ex}{1ex plus .2ex}{\bf}}

\newcommand{\amgrp}[1]{{\mathbf{#1}}}
\newcommand{\comp}[1]{{\mathrm{#1}}}

\newcommand{\amgrpA}{{\mathbf{A}}}
\newcommand{\amgrpC}{{\mathbf{C}}}
\newcommand{\amgrpD}{{\mathbf{D}}}
\newcommand{\amgrpB}{{\mathbf{B}}}

\newcommand{\amgrpG}{{\mathbf{G}}}
\newcommand{\amgrpH}{{\mathbf{H}}}

\newcommand{\amgrpP}{{\mathbf{P}}}
\newcommand{\amgrpT}{{\mathbf{T}}}
\newcommand{\amgrpU}{{\mathbf{U}}}
\newcommand{\amgrpX}{{\mathbf{X}}}

\newcommand{\ama}{{\mathbf a}}

\newcommand{\amg}{{\mathbf g}}
\newcommand{\famg}{\ul{\mathbf g}}
\newcommand{\amh}{{\mathbf h}}

\newcommand{\bamgrpA}{\bar{\amgrpA}}

\newcommand{\bamgrpG}{\bar{\amgrpG}}
\newcommand{\bamgrpH}{\bar{\amgrpH}}
\newcommand{\bamgrpX}{\bar{\amgrpX}}

\newcommand{\compA}{{{A}}}

\newcommand{\compG}{{{G}}}
\newcommand{\compH}{{{H}}}

\newcommand{\compX}{{{X}}}
\newcommand{\compa}{\alpha}

\newcommand{\compg}{\gamma}

\newcommand{\ucompA}{{\tilde{A}}}

\newcommand{\ucompa}{{\tilde{\alpha}}}

\newcommand{\umap}{\pi}

\def\<{\langle}
\def\>{\rangle}

\DeclareMathOperator{\Tr}{Tr}

\newcommand{\amA}{{\mathscr{A}}}

\newcommand{\famG}{\ul{\mathscr{G}}}
\newcommand{\amG}{{\mathscr{G}}}
\newcommand{\amH}{{\mathscr{H}}}

\DeclareMathOperator{\Spin}{Spin}
\usepackage[refpage]{nomencl}
\makenomenclature		
\newcommand{\nom}{}

\DeclareMathOperator{\edg}{E}
\DeclareMathOperator{\vrtc}{V}
\newcommand{\coxdiag}{\Delta}
\newcommand{\liediag}{\Gamma}

\begin{document}
\title[Curtis-Tits and Phan amalgams]{Classification of Curtis-Tits and Phan amalgams with $3$-spherical diagram}
\author[Blok]{Rieuwert J. Blok}
\address{Department of Mathematics and Statistics\\
Bowling Green State University\\
Bowling Green, oh 43403\\
U.S.A.
}
\curraddr{School of Mathematics\\
University of Birmingham\\
Edgbaston, B15 2TT\\
U.K.
}
\email{blokr@member.ams.org}

\author[Hoffman]{Corneliu G. Hoffman}
\address{School of Mathematics\\
University of Birmingham\\
Edgbaston, B15 2TT\\
U.K.
}
\email{C.G.Hoffman@bham.ac.uk}

\author[Shpectorov]{Sergey V. Shpectorov}
\address{School of Mathematics\\
University of Birmingham\\
Edgbaston, B15 2TT\\
U.K.}
\email{S.V.Shpectorov@bham.ac.uk}

\begin{abstract}
We classify all non-collapsing Curtis-Tits and Phan amalgams with $3$-spherical diagram over all fields.
In particular, we show that amalgams with spherical diagram are unique, a result required by the classification of finite simple groups. 
We give a simple  condition on the amalgam which is necessary and sufficient for it to arise from a group of Kac-Moody type. This also yields a definition of a large class of groups of Kac-Moody type in terms of a finite presentation.
\end{abstract}
\maketitle

%\begin{keyword}
% Curtis-Tits amalgam, Phan amalgam, classification and uniqueness, group of Lie type.
%%    AMS subject classification (2000):
%\MSC{2010} 20G35 \sep % Linear algebraic groups over ad`eles and other rings and schemes
%%    Secondary 
%%51E24%
%%    %20G15, %linear algebraic groups over arbitrary fields
%%    %20G40,% linear algebraic groups over finite fields
%%    %20E42,%Groups with a $BN$-pair; buildings (as method for classification)
%%    %51E24, % buildings and the geometry of diagrams
%%    %Topological methods in group theory
%\end{keyword}%\date{}							% Activate to display a given date or no date

%\tableofcontents
%\section{}
%\subsection{}
\section{Introduction}
Local recognition results play an important role in various parts of mathematics. A key example comes from the monumental classification of finite simple groups. Local analysis of the unknown finite simple group $G$ yields a local datum consisting of a small collection of subgroups fitting together in a particular way, called an amalgam. 
The Curtis-Tits theorem~\cite{Cur1965a, Tim1998,Tim03,Tim04,Tim06} and the Phan (-type) theorems~\cite{Pha1971,Pha1977,Pha1977a} describe amalgams appearing in known groups of Lie type.
Once the amalgam in $G$ is identified as one of the amalgams given by these theorems, $G$ is known. 

The present paper was partly motivated by a question posed by R.~Solomon and R.~Lyons about this identification step, arising from their work on the classification~\cite{Gor1983,GorLyoSol1996,GorLyoSol1998,GorLyoSol1999,GorLyoSol2002,GorLyoSol2005}: 
%The present paper was partly motivated by a question posed by R.~Solomon and R.~Lyons about this identification step: 
Are Curtis-Tits and Phan type amalgams uniquely determined by their subgroups? More precisely is there a way of fitting these subgroups together so that the amalgam gives rise to a different group? 
In many cases it is known that, indeed, depending on how one fits the subgroups together, either the resulting amalgam arises from these theorems, or it does not occur in any non-trivial group.
This is due to various results of Bennet and Shpectorov~\cite{BeSh2004}, Gramlich~\cite{Gr2004}, Dunlap~\cite{Dun2005}, and R.~Gramlich, M.~Horn, and W.~Nickel~\cite{GrHoNi2006}.
However, all of these results use, in essence, a crucial observation by Bennett and Shpectorov about tori in rank-$3$ groups of Lie type, which fails to hold for small fields.
In the present paper we replace the condition on tori by a more effective condition on root subgroups, which holds for all fields. This condition is obtained by a careful analysis of maximal subgroups of groups of Lie type.
Thus the identification step can now be made for all possible fields.
A useful consequence of the identification of the group $G$, together with the Curtis-Tits and Phan type theorems, is that it yields a simplified version of the Steinberg presentation for $G$.

Note that this solves the  - generally much harder - existence problem: ``how can we tell if a given amalgam appears in any non-trivial group?''

\medskip
The unified approach in the present paper not only extends the various results on Curtis-Tits and Phan amalgams occurring in groups of Lie type to arbitrary fields, but in fact also applies to a much larger class of Curtis-Tits and Phan type amalgams, similar to those occurring in groups of Kac-Moody type. Here, both the uniqueness and the existence problem become significantly more involved.

Groups of Kac-Moody type were introduced by J.~Tits as automorphism groups of certain geometric objects called twin-buildings \cite{Ti1992}. In the same paper J.~Tits conjectured that these groups are uniquely determined by the fact that the group acts on some twin-building, together with local geometric data called Moufang foundations. As an example he sketched an approach towards classifying such foundations in the case of simply-laced diagrams.
This conjecture was subsequently proved for Moufang foundations built from locally split and locally finite rank-$2$ residues by B.~M\"uhlherr in~\cite{Mu1999} and refined by P.~E.~Caprace in~\cite{Cap2007}.
All these results produce a classification of groups of Kac-Moody type using local data in the form of an amalgam, together with a global geometric assumption stipulating the existence of a twin-building on which the group acts.

Ideally, one would use the generalizations of the Curtis-Tits and Phan type theorems to describe the groups of Kac-Moody type in terms of a simplified Steinberg type presentation. 
However, the geometric assumption is unsatisfactory for this purpose as it is impossible to verify directly from the presentation itself. 

In our unified approach we consider all possible amalgams whose local structure is any one of those appearing in the above problems. There is no condition on the field.
Then, we classify those amalgams that satisfy our condition on root groups and show that in the spherical case they are unique. This explains why groups of Lie type can uniquely be recognized by their amalgam.
By contrast, in the non-spherical case the amalgams are not necessarily unique and, indeed, not all such amalgams give rise to groups of Kac-Moody type. This is a consequence of the fact that we impose no global geometric condition.
Nevertheless, we give a simple condition on the amalgam itself which decides whether it comes from a group of Kac-Moody type or not. As a result, we obtain a purely group theoretic definition of a large class of groups of Kac-Moody type just in terms of a finite presentation.

Finally, we note that an amalgam must satisfy the root subgroup condition to occur in a non-trivial group.
A subsequent study generalizing~\cite{BloHof2014a,BloHof2016} shows that in fact all amalgams satisfying the root group condition do occur in non-trivial groups. Thus, in this much wider context the existence problem is also solved.

\medskip
We shall now give an overview of the results in the present paper.
Recall that a Dynkin diagram $\liediag$ is an oriented edge-labelled graph.
We say that $\liediag$ is {\dfn connected} if the underlying (unlabelled) graph is connected in the usual sense.
Moreover, we use topological notions such as spanning tree and homotopy rank of $\liediag$ referring to the underlying graph.

For Phan amalgams we prove the following (for the precise statement see  Theorem~\ref{thm:P classification of 3-spherical amalgams}).

\bmth\label{mthm:P classification of 3-spherical amalgams}
Let  $q$ be any prime power and let $\liediag$ be a connected $3$-spherical diagram with homotopy rank $r$.
Then, there is a bijection between the elements of $\prod_{s=1}^r \Aut(\FF_{q^2})$ and the type preserving isomorphism classes of  Curtis-Tits amalgams with diagram $\liediag$ over $\FF_q$.
\emth

For Curtis-Tits amalgams the situation is slightly more complicated (for the precise statement see Theorem~\ref{thm:CT classification of 3-spherical amalgams}).

\bmth\label{mthm:CT classification of 3-spherical amalgams}
Let  $q$ be a prime power and let $\liediag$ be a connected $3$-spherical diagram with homotopy rank $r$.
Then there exists a set of positive integers $\{e_1,\ldots,e_r\}$ so that 
there is a bijection between the elements of $\prod_{s=1}^r \Aut(\FF_{q^{e_s}})\times \ZZ/2\ZZ$ and the type preserving isomorphism classes of  Curtis-Tits amalgams with diagram $\liediag$ over $\FF_q$.
\emth

\bmco\label{maincor:tree diagrams}
Let $q$ be a prime power and let $\liediag$ be a $3$-spherical tree. Then, up to type preserving isomorphism, there is a unique Curtis-Tits and a unique Phan amalgam over $\FF_q$ with diagram $\liediag$.
\emco

Note that Corollary~\ref{maincor:tree diagrams} includes all spherical diagrams of rank $\ge 3$.
Several special cases of the above results were proved elsewhere. Indeed, Theorem~\ref{mthm:CT classification of 3-spherical amalgams} was proved for simply-laced diagrams and $q\ge 4$ in~\cite{BloHof2014b}.
 Corollary~\ref{maincor:tree diagrams}  was proved for Phan amalgams with $\liediag=A_n$ in~\cite{BeSh2004}, for general simply-laced tree diagram in~\cite{Dun2005}, and for $\liediag=C_n$ for $q\ge 3$ in~\cite{Gra2004,GraHorNic2006}.

\medskip
The classification of Curtis-Tits amalgams will be done along the following lines. Note if $(\amgrpG,\bamgrpG_i,\bamgrpG_j)$ is a Curtis-Tits standard of type different from $A_1\times A_1$, and $\bamgrpX$ is any Sylow $p$-subgroup in 
one of the vertex groups, say $\bamgrpG_i$, then generically it generates $\amgrpG$ together with $\bamgrpG_j$.
In Subsection~\ref{subsub:fundamental root groups in CT standard pairs} we show that there is a unique pair $(\bamgrpX_i^+,\bamgrpX_i^-)$ of Sylow $p$-subgroups in $\bamgrpG_i$ whose members do not have this property. Moreover, each member commutes with a unique member in the other vertex group.

In Subsection~\ref{subsub:weak systems} we show that in a non-collapsing Curtis-Tits amalgam $\amG=\{\amgrpG_{i},\amgrpG_{i, j}, \amg_{i,j} \mid  i, j \in I\}$ for each $i$ there exists a pair $(\amgrpX_i^+,\amgrpX_i^-)$ of Sylow subgroups in $\amgrpG_i$ such that for any edge $\{i,j\}$ $(\amg_{i,j}(\amgrpX_i^+, \amgrpX_i^-)$ is the pair for $(\amgrpG_{i,j},\bamgrpG_i,\bamgrpG_j)$ as above. The collection $\Chi=\{\amgrpX_i^+,\amgrpX_i^-\colon i\in I\}$ is called a {\dfn weak system of fundamental root groups}. Without loss of generality one can assume that any amalgam with the same diagram has the exact same weak system $\Chi$.
As a consequence all amalgams with the same diagram can be determined up to isomorphism by studying the coefficient system associated to $\Chi$, that is, the graph of groups consisting of automorphisms of the vertex and edge groups preserving $\Chi$.
In Subsection~\ref{subsec:coefficient system} we determine the coefficient system associated to $\Chi$.
In Subsection~\ref{subsec:trivial support}, we pick a spanning tree $\Sigma$ for $\liediag$ and use precise information about the coefficient system to create a standard form of a Curtis-Tits amalgam in which all vertex-edge inclusion maps are trivial except for the edges in $\Sigma$.
In particular this shows that if $\liediag$ is a tree, then the amalgam is unique up to isomorphism.
Finally in Subsection~\ref{subsec:classification of 3-spherical CT amalgams} we show that for a suitable choice of $\Sigma$, the remaining non-trivial inclusion maps uniquely determine the amalgam.

The classification of Phan amalgams  in Section~\ref{sec:classification of Ph amalgams} follows the same pattern. However in this case the role of the weak system of fundamental root groups is replaced by a system of tori in the vertex groups, whose images in the edge groups must form a torus there.

\medskip
As shown here, the existence of a weak system of fundamental root groups is a necessary condition for the existence of a non-trivial completion. A natural question of course is whether it is also sufficient.
In the spherical cases, the amalgams are unique and the Curtis-Tits and Phan theorems identify universal completions of these amalgams.
In~\cite{BloHof2016} it is shown that any Curtis-Tits amalgam with $3$-spherical simply-laced diagram over a field with at least four elements having property (D) has a non-trivial universal completion, which is identified up to a rather precisely described central extension.
In the present paper we will not study completions of the Curtis-Tits and Phan amalgams classified here, but merely note that similar arguments yield non-trivial completions for all amalgams. In particular, the conditions mentioned above are indeed sufficient for the existence of these completions.
In general we don't know of a direct way of giving conditions on an amalgam ensuring the existence of a non-trivial completion.

\section*{Acknowledgements}
This paper was written as part of the project KaMCAM funded by the European Research Agency through a Horizon 2020 Marie-Sk\l odowska Curie fellowship (proposal number 661035).

\section{Curtis-Tits and Phan amalgams and their diagrams}
\subsection{Diagrams}\label{subsec:diagrams}
In order to fix some notation, we start with some definitions.

\bde\label{dfn:coxdiag}
A {\dfn Coxeter matrix  over the set $I=\{1,2,\ldots,n\}$} of finite cardinality $n$ is a symmetric matrix 
$M=(m_{ij})_{i,j\in I}$ with entries in $\NN_{\ge 1}\cup \{\infty\}$ such that, for all $i,j\in I$ distinct we have $m_{ii}=1$ and $m_{ij}\ge 2$.

A {\dfn Coxeter diagram with Coxeter matrix $M$} is an edge-labelled graph $\coxdiag=(I,E)$ with vertex set $I=\vrtc \coxdiag$ and edge-set $E=\edg\coxdiag$ without loops such that 
 for any distinct $i,j\in I$, there is an edge labelled $m_{ij}$ between $i$ and $j$ whenever $m_{ij}>2$; if $m_{i,j}=2$, there is no such edge.
Thus, $M$ and $\coxdiag$ determine each other uniquely.
For any subset $J\sbe I$, we let $\coxdiag_J$ denote the diagram induced on vertex set $J$.
We say that $\coxdiag$ is {\dfn connected} if the underlying (unlabelled) graph is connected in the usual sense.
Moreover, we use topological notions such as spanning tree and homotopy rank of $\coxdiag$ referring to the underlying graph.

A {\dfn Coxeter system with Coxeter matrix $M$} is a pair $(W, S)$, where $W$ is a group generated by the set $S=\{s_i\colon i\in I\}$ subject to the relations $(s_i s_j)^{m_{ij}}=1$ for all $i,j\in I$.
For each subset $J\sbe I$, we let $W_J=\langle s_j\colon j\in J\rangle_W$.
We call $M$ and $(W,S)$ $m$-spherical if every subgroup $W_J$ with $|J|=m$ is finite ($m\in \NN_{\ge 2}$).
Call $(W,S)$ spherical if it is $n$-spherical.
\ede

In order to describe Curtis-Tits and Phan amalgams more precisely, we also introduce a Lie diagram.
\bde\label{dfn:liediag}
Let $\coxdiag=(I,E)$ be a Coxeter diagram. 
A {\dfn Lie diagram of Coxeter type $\coxdiag$} is an untwisted or twisted Dynkin diagram $\liediag$ whose edge labels $l_{ij}$  do not specify the orientation.
In this paper we shall only be concerned with Lie diagrams of Coxeter type $A_n$, $B_n$, $D_n$, $E_6$, $E_7$, $E_8$, and $F_4$.
For these, we have the following correspondence   
$$
{\footnotesize
\begin{array}{c||c|c|c|c|c}
\coxdiag &  A_n & B_n  & D_n, & E_n (n=6,7,8) & F_4   \\
\hline
\liediag  & A_n & B_n, C_n,  \twD_{n+1}, \twA_{2n-1}, \twA_{2n} & D_n & E_n  & F_4, F_4^*,  \twE_6 ,\twE_6^*
\end{array}
}$$
Here $F_4$ and $\twE_6$ (resp. $F_4^*$ and $\twE_6^*$) denote the diagrams where node $1$ corresponds to the long (resp.~short) root (Bourbaki labeling).
\ede

Let us introduce some more notation.
We shall denote the Frobenius automorphism of order $2$ of $\FF_{q^2}$ by $\sigma$.
Below we will consider sesquilinear forms $\hform$ on an $\FF_{q^2}$-vector space $V$.
By convention, all these forms are linear in the first coordinate, that is
 $\hform(\lambda u,\mu v)=\lambda\hform(u,v)\mu^\sigma$ for $u,v\in V$ and $\lambda, \mu\in \FF_{q^2}$. 
Recall that $\hform$ is hermitian if $\hform(v,u)=\hform(u,v)^\sigma$ for all $u,v\in V$.

\subsection{Standard pairs of Curtis-Tits type}\label{subsec:standard CT pairs}
Let $\liediag$ be a Lie diagram of type $A_2$, $B_2/C_2$, $\twD_3/\twA_3$ and $q=p^e$ for some prime $p\in \ZZ$ and $e\in \ZZ_{\ge 1}$. Then a  {\dfn Curtis-Tits standard pair of type $\liediag(q)$}
 is a triple $(\amgrpG,\amgrpG_1,\amgrpG_2)$ of groups such that one of the following occurs:
\paragraph{\bf ($\liediag=A_1\times A_1$)} 
Now $\amgrpG=\amgrpG_1\times \amgrpG_2$ and $\amgrpG_1\cong\amgrpG_2\cong\SL_2(q)$.

\paragraph{\bf($\liediag=A_2$)} 
Now $\amgrpG=\SL_3(q)=\SL(V)$ for some $\FF_q$-vector space with basis $\{e_1,e_2,e_3\}$ and  
 $\amgrpG_1$ (resp. $\amgrpG_2$) is the stabilizer of the subspace $\langle e_1,e_2\rangle$ (resp. $\langle e_2,e_3\rangle$) and the vector $ e_3$ (resp. $ e_1$).
 
Explicitly we have
 \begin{align*}
\amgrpG_1&=\left\{\begin{pmatrix} 
a & b &  \\
c & d   & \\
 &  & 1 \\
\end{pmatrix}\colon a,b,c,d\in\FF_{q} \mbox{ with } ad-bc=1\right\}, \\
\amgrpG_2&=\left\{\begin{pmatrix} 
1 &   & \\
 & a  & b \\
 &  c & d  \\
\end{pmatrix}\colon a,b,c,d\in\FF_{q} \mbox{ with } ad-bc=1  \right\}.
  \end{align*}

\paragraph{\bf($\liediag=C_2$)}Now $\amgrpG=\Sp_4(q)=\Sp(V,\beta)$, where $V$ is an $\FF_q$-vector space with basis $\{e_1,e_2,e_3,e_4\}$ and 
$\beta$ is the symplectic form with Gram matrix
 \begin{align*}
 M & = \begin{pmatrix} 
  & & 1 & 0 \\ 
  & & 0 & 1 \\ 
  -1 & 0&  &  \\ 
  0 & -1 &  &  
\end{pmatrix}.
 \end{align*}
 $\amgrpG_1\cong \SL_2(q)$ is the derived subgroup of $\Stab_\amgrpG(\langle e_1,e_2\rangle)\cap \Stab_\amgrpG(\langle e_3,e_4\rangle)$ and $\amgrpG_2=\Stab_\amgrpG(e_1)\cap \Stab_\amgrpG(e_3)\cong \Sp_2(q)\cong\SL_2(q)$.
Explicitly we have
 \begin{align*}
\amgrpG_1&=\left\{\begin{pmatrix} 
a & b & & \\
c & d  & & \\
 &  & d & -c \\
 &  &-b & a \end{pmatrix}\colon a,b,c,d\in\FF_{q} \mbox{ with } ad-bc=1  \right\},\\
\amgrpG_2&=\left\{\begin{pmatrix} 
1 &  & & \\
 & a  & & b\\
 &  & 1& \\
 & c & &d \end{pmatrix}\colon a,b,c,d\in\FF_{q} \mbox{ with } ad-bc=1\right\}.
  \end{align*}

\bre\label{rem:B2 q odd}
We are only interested in Curtis-Tits standard pairs of type $B_2$ for $q$ odd. However, in that case we have $\Spin_{5}(q)\cong\Sp_4(q)$ is the unique central extension of the simple group $\Omega_5(q)\cong \PSp_4(q)$. Therefore, we can also describe the Curtis-Tits standard pair for $B_2$ as a Curtis-Tits standard pair for $C_2$ with $\amgrpG_1$ and $\amgrpG_2$ interchanged.
\ere
\paragraph{\bf($\liediag=\twA_3$)} Now $\amgrpG=\SU_4(q)=\SU(V)$ for some $\FF_{q^2}$-vector space $V$ with basis $\{e_1,e_2,e_3,e_4\}$ equipped with a non-degenerate hermitian form 
 $\hform$ for which this basis is hyperbolic with Gram matrix
 \begin{align*}
 M&=\begin{pmatrix}
  & & 1 & 0 \\
  & & 0 &1 \\
  1 & 0 &  &  \\
 0 & 1 &  &  
 \end{pmatrix}.
 \end{align*}
Now  $\amgrpG_1$ is the derived subgroup of the simultaneous stabilizer of the subspaces $\langle e_1,e_2\rangle$ and $\langle e_3,e_4\rangle$ and $\amgrpG_2$ is the stabilizer of the vectors $e_1$ and $e_3$ and the hyperbolic line $\langle e_2,e_4\rangle$.
We have $\amgrpG_2\cong\SU_2(q)\cong \SL_2(q)$ and $\amgrpG_1\cong \SL_2(q^2)$.
Explicitly we have
 \begin{align*}
\amgrpG_1&=\left\{\begin{pmatrix} 
a & b & & \\
c & d  & & \\
 &  & d^\sigma & -c^\sigma \\
 &  &-b^\sigma & a^\sigma \end{pmatrix}\colon a,b,c,d\in\FF_{q^2} \mbox{ with } ad-bc=1  \right\}, \\
\amgrpG_2&=\left\{\begin{pmatrix} 
1 &  & & \\
 & a  & & b\eta\\
 &  & 1& \\
 & c\eta^{-1} & &d \end{pmatrix}\colon a,b,c,d\in\FF_{q} \mbox{ with } ad-bc=1\right\},\\
   \end{align*}
 where $\eta\in \FF_{q^2}$   has $\eta+\eta^q=0.$

\bre\label{rem:twD3}
For completeness we also define a standard Curtis-Tits pair $(\amgrpH,\amgrpH_1,\amgrpH_2)$ of type $\twD_3(q)$.
Take $\amgrpH=\Omega_6^-(q)=\Omega^-(V,\cQ)$,  where $V$ is an $\FF_q$-vector space with basis $\{e_1,$ $e_2,$ $e_3,$ $e_4,$ $e_5,$ $e_6\}$ and 
$\cQ(\sum_{i=1}^5 x_i e_i)=x_1x_3+x_2x_4+f(x_5,x_6)$, for some quadratic polynomial $f(x,1)$ that is irreducible over $\FF_q$.
Here, 
 $\amgrpH_1\cong \SL_2(q)$ is the derived subgroup of $\Stab_\amgrpH(\langle e_1,e_2\rangle)\cap \Stab_\amgrpH(\langle e_3,e_4\rangle)$ if $\SL_2(q)$ is perfect that is $q>2$, and it is the subgroup $\Stab_\amgrpH(\langle e_1,e_2\rangle)\cap \Stab_\amgrpH(\langle e_3,e_4\rangle)\cap \Stab(v)$ for some non-singular vector $v\in \langle e_5,e_6\rangle$ if $q=2$, and $\amgrpH_2=\Stab_\amgrpH(e_1)\cap \Stab_\amgrpH(e_3)\cong \Omega_4^-(q)\cong\PSL_2(q^2)$.

However, there exists a unique standard Curtis-Tits pair $(\amgrpG,\amgrpG_1,\amgrpG_2)$ of type $\twA_3(q)$ and a surjective homomorphism $\pi\colon \amgrpG\to \amgrpH$ with $\ker\pi=\{\pm 1\}=Z(\amgrpG_1)\le Z(\amgrpG)$. 
It induces $\pi\colon \amgrpG_1\cong\SL_2(q^2)\to \Omega_4^-(q)=\amgrpH_2\cong \PSL_2(q^2)$ and 
$\pi\colon\amgrpG_2\cong \SU_2(q)\to \SL_2(q)=\amgrpH_2$.
Because of this map, any amalgam involving a standard Curtis-Tits pair of type $\twD_3(q)$ is the image of an amalgam involing a standard Curtis-Tits pair of type $\twA_3(q)$ (see~Subsection~\ref{subsec:CTP amalgams}).
\ere

\bde\label{dfn:standard CT identification map}
For Curtis-Tits amalgams, the {\dfn standard identification map} will be the isomorphism $\amg\colon \SL_2(q^e) \to \amgrpG_i$ sending 
\begin{align*}
\begin{pmatrix}
a & b \\
c & d\end{pmatrix}
\end{align*} 
to the corresponding matrix of $\amgrpG_i$ as described above.
Here $e=1$ unless $\liediag(q)=\twA_3(q)$ and $i=1$ or $\liediag(q)=\twD_3$ and $i=2$, in which case $e=2$.
\ede

\subsection{Standard pairs of Phan type}\label{subsec:standard P pairs}

Let $\liediag$ be as above. Then a  {\dfn  Phan standard pair of type $\liediag(q)$}
 is a triple $(\amgrpG,\amgrpG_1,\amgrpG_2)$ such that one of the following occurs:

\paragraph{\bf ($\liediag=A_1\times A_1$)}
Now $\amgrpG=\amgrpG_1\times \amgrpG_2$ and $\amgrpG_1\cong\amgrpG_2\cong\SU_2(q)=\SU(V)$ for some $\FF_{q^2}$-vector space $V$ with basis $\{e_1,e_2\}$ equipped with a non-degenerate hermitian form 
 $\hform$ for which this basis is orthonormal.
\paragraph{\bf($\liediag=A_2$)}Now $\amgrpG=\SU_3(q)=\SU(V)$ for some $\FF_{q^2}$-vector space $V$ with basis $\{e_1,e_2,e_3\}$ equipped with a non-degenerate hermitian form 
 $\hform$ for which this basis is orthonormal.
As in the Curtis-Tits case, $\amgrpG_1$ (resp. $\amgrpG_2$) is the stabilizer of the subspace  $\langle e_1,e_2\rangle$ (resp. $\langle e_2,e_3\rangle$) and the vector $ e_3$ (resp. $ e_1$).
 We have $\amgrpG_1\cong\amgrpG_2\cong \SU_2(q)$.
 
 Explicitly we have
 \begin{align*}
\amgrpG_1&=\left\{\begin{pmatrix} 
a & b & & \\
-b^\sigma & a^\sigma   & \\
 &   & 1 \end{pmatrix}\colon a,b\in\FF_{q^2} \mbox{ with } aa^\sigma+bb^\sigma=1  \right\}, \\
\amgrpG_2&=\left\{\begin{pmatrix} 
1 &  &  \\
 & a   & b\\
 & -b^\sigma &a^\sigma\end{pmatrix}\colon a,b\in\FF_{q^2} \mbox{ with } aa^\sigma+bb^\sigma=1\right\}. \end{align*}

\paragraph{\bf($\liediag=C_2$)}  Let $V$ be an $\FF_{q^2}$-vector space with basis $\{e_1,$ $e_2,$ $e_3,$ $ e_4\}$ and let 
$\beta$ be the symplectic form with Gram matrix
 \begin{align*}
 M & = \begin{pmatrix} 
  & & 1 & 0 \\ 
  & & 0 & 1 \\ 
  -1 & 0&  &  \\ 
  0 & -1 &  &  
\end{pmatrix}.
 \end{align*}
Moreover, let $\hform$ be the (non-degenerate) hermitian form 
 $\hform$ for which this basis is orthonormal.
 
 Now $\amgrpG=\Sp(V,\beta)\cap \SU(V,\hform)\cong \Sp_4(q)$, $\amgrpG_1\cong \SU_2(q)$ is the derived subgroup of $\Stab_\amgrpG(\langle e_1,e_2\rangle)\cap \Stab_\amgrpG(\langle e_3,e_4\rangle)$ and $\amgrpG_2=\Stab_\amgrpG(e_1)\cap \Stab_\amgrpG(e_3)\cong \Sp_2(q)\cong\SU_2(q)$, .
Note that $Z(\amgrpG) =Z(\amgrpG_1)$ and $Z(\amgrpG)\cap \amgrpG_2=\{1\}$.

Explicitly we have
 \begin{align*}
\amgrpG_1&=\left\{\begin{pmatrix} 
a & b & & \\
-b^\sigma & a^\sigma  & & \\
 &  & a^\sigma & b^\sigma \\
 &  &-b &a \end{pmatrix}\colon a,b\in\FF_{q^2} \mbox{ with } aa^\sigma+bb^\sigma=1  \right\},\\
\amgrpG_2&=\left\{\begin{pmatrix} 
1 &  & & \\
 & a  & & b\\
 &  & 1& \\
 & -b^\sigma & &a^\sigma\end{pmatrix}\colon a,b\in\FF_{q^2} \mbox{ with } aa^\sigma+bb^\sigma=1\right\}.   \end{align*}

\bde\label{dfn:standard Phan identification map}
For Phan amalgams, the {\dfn standard identification map} will be the isomorphism $\amg\colon \SU_2(q) \to \amgrpG_i$ sending 
\begin{align*}
\begin{pmatrix}
a & b \\
-b^\sigma & a^\sigma\end{pmatrix}
\end{align*} 
to the corresponding matrix of $\amgrpG_i$ as described above.
\ede

\subsection{Amalgams of Curtis-Tits and Phan type}\label{subsec:CTP amalgams}
\bde\label{dfn:amalgam}
An \nom{{\em amalgam} }{} over a poset $(\cP,\prec)$ is a collection $\amA=\{\amgrpA_x\mid x\in \cP\}$ of groups, together with a collection $\ama_\bullet=\{\ama_x^y\mid x\prec y, x,y\in \cP\}$ of monomorphisms $\ama_x^y\colon \amgrpA_x\into \amgrpA_y$, called {\em inclusion maps} such that whenever $x\prec y\prec z$, we have
 $\ama_x^z=\ama_y^z\after\ama_x^y$; we shall write $\bamgrpA_x=\ama_x^y(\amgrpA_x)\le \amgrpA_y$. 
A {\em completion} of $\amA$ is a group $\compA$ together with a collection  $\compa_\bullet=\{\compa_x\mid x\in \cP\}$ of homomorphisms $\compa_x\colon \amgrpA_x\to \compA$, whose images - often denoted ${\compA}_x=\compa_x(\amgrpA_x)$ - generate $\compA$, such that for any $x, y\in \cP$ with $x\prec y$ we have 
$\compa_y\after\compa_x^y=\compa_x$.
The amalgam $\amA$ is {\em non-collapsing} if it has a non-trivial completion.
As a convention, for any subgroup $\amgrp{H}\le \amgrpA_J$, let $\comp{H}=\alpha(\amgrpH)\le \compA$.

A completion $(\ucompA,\ucompa_\bullet)$ is called {\em universal} if for any completion $(\compA,\compa_\bullet)$ there is a unique surjective group homomorphism $\umap\colon \ucompA\to \compA$ such that $\compa_\bullet=\umap\after\ucompa_\bullet$. A universal completion always exists.

\ede
\bde\label{dfn: CTP structure}
Let $\liediag=(I, E)$ be  a Lie  diagram.
A  {\em   Curtis-Tits (resp.~Phan)  amalgam with Lie diagram  $\liediag$ over $\FF_q$}  is an amalgam 
$\amG=\{\amgrpG_{i},\amgrpG_{i, j}, \amg_{i,j} \mid  i, j \in I\}$ 
over $\cP=\{J\mid \emptyset\ne J\sbe I \mbox{ with }|J|\le 2\}$ ordered by inclusion such that 
for every $i,j\in I$,
$(\amgrpG_{i,j}, \bamgrpG_i, \bamgrpG_j)$ is a Curtis-Tits / Phan  standard pair of type $\liediag_{i,j}(q^e)$, for some $e\ge 1$ as defined in Subsection~\ref{subsec:standard CT pairs}~and~\ref{subsec:standard P pairs}.
Moreover $e=1$ is realized for some $i,j\in I$.
Note that in fact $e$ is always a power of $2$.
This follows immediately from  connectedness of the diagram and the definition of the standard pairs of type $A_2$, $C_2$, and $\twA_3$.
For any subset $K\sbe I$ , we let 
\begin{align*}
\amG_K&=\{\amgrpG_{i},\amgrpG_{i, j}, \amg_{i,j} \mid  i, j \in K\}. 
\end{align*}
 \ede
 
\bre 
Suppose that one considers an amalgam 
\begin{align*}
\amH=\{\amgrpH_{i},\amgrpH_{i, j}, \amh_{i,j} \mid  i, j \in I\}
\end{align*} 
  over $\FF_q$ with diagram $\liediag$, such that for any $i,j\in I$, the triple $(\amgrpH_{i,j},$ $\bamgrpH_i, $ $\bamgrpH_j)$ 
  is not a standard pair, but there is a standard pair $(\amgrpG_{i,j},\bamgrpG_i,\bamgrpG_j)$ such that the respective $\amgrpH$'s are central quotients of the corresponding $\amgrpG$'s.
Then, $\amH$ is the quotient of a unique Curtis-Tits or Phan amalgam over $\FF_q$ with diagram $\liediag$.
 Hence for classification purposes it suffices to consider Curtis-Tits or Phan amalgams. 
 In particular, in view of Remark~\ref{rem:twD3}, we can restrict ourselves to Curtis-Tits amalgams in which the only rank-$2$ subdiagrams are of type $A_1\times A_1$, $A_2$, $C_2$, and $\twA_3$.
\ere

\bde\label{dfn:isomorphism of amalgams}
Suppose $\amG=\{\amgrpG_{i},\amgrpG_{i, j}, \amg_{i,j} \mid  i, j \in I\}$, $\amG^+=\{\amgrpG_{i}^+,\amgrpG_{i, j}^+, \amg_{i,j}^+ \mid  i, j \in I\}$  are two Curtis-Tits (or Phan) amalgams over $\FF_q$ with the same diagram $\liediag$.
Then a type preserving isomorphism $\phi\colon \amG\to\amG^+$ is a collection $\phi=\{\phi_i,\phi_{i,j}\colon i,j\in I\}$ of group isomorphisms such that, for all $i,j\in I$, we have  
\begin{align*}
\phi_{i,j}\after \amg_{i,j} & = \amg_{i,j}^+\after \phi_i\\
\phi_{i,j}\after \amg_{j,i} & = \amg_{j,i}^+\after \phi_j.
\end{align*}
Unless indicated otherwise, this is the kind of isomorphism we shall consider, omitting the term "type preserving". 
It is also possible to consider type permuting isomorphisms, defined in the obvious way.
\ede

\section{Background on groups of Lie type}
\subsection{Automorphisms of groups of Lie type of small rank}\label{subsec:automorphism groups}
Automorphisms of groups of Lie type are all known. In this subsection we collect some facts that we will need later on.
We shall use the notation from~\cite{Wil2009}.
\paragraph{Automorphisms of $\SL_n(q)$}
Define automorphisms of $\SL_n(q)$ as follows (where $x=(x_{ij})_{i,j=1}^n\in \SL_n(q)$):
\begin{align*}
c_g& \colon x\mapsto x^g=g^{-1}xg && (g\in \PGL_n(q)),\\
\alpha & \colon x\mapsto x^\alpha=(x_{ij}^\alpha)_{i,j=1}^n && (\alpha\in \Aut(\FF_q)),\\
\trin&\colon x\mapsto x^\trin = {}^t x^{-1} &&(\mbox{transpose-inverse}).
\end{align*}
We note that for $n=2$, $\trin$ coincides with the map 
 $x\mapsto x^\mu$, where
$\mu =\begin{pmatrix}  
 0 & -1 \\
 1 & 0\end{pmatrix}$.
We let $\PGammaL_n(q)=\PGL_n(q)\rtimes \Aut(\FF_q)$.
 
\paragraph{Automorphisms of $\Sp_{2n}(q)$}
Outer automorphisms of $\Sp_{2n}(q)$ are of the form $\Aut(\FF_q)$ as for $\SL_{2n}(q)$, defined with respect to a symplectic basis, or come from the group $\GSp_{2n}(q)\cong \Sp_{2n}(q)\mspl (\FF_q^*/(\FF_q^*)^2)$ of linear similarities of the symplectic form, where 
  $\FF_q^*$ acts as conjugation by 
  \begin{align*}
  \delta(\lambda) &= \begin{pmatrix}
  \lambda I_n & 0_n\\
  0_n & I_n
  \end{pmatrix}
&& (\lambda\in \FF_q^*)
  \end{align*}
This only provides a true outer automorphism if $\lambda$ is not a square and we find that 
 $\PGSp_{2n}(q)\cong \PSp_{2n}(q)\mspl 2$ if $q$ is odd and $\PGSp_{2n}(q)=\PSp_{2n}(q)$ if $q$ is even.
We define 
\begin{align*}
\GammaSp_{2n}(q)&= \GSp_{2n}(q)\rtimes \Aut(\FF_q)\\
\PGammaSp_{2n}(q)&=\PGSp_{2n}(q)\rtimes \Aut(\FF_q).
\end{align*}
Note that, as in $\SL_2(q)$, the map $\trin\colon A\to {}^tA^{-1}$ is the inner automorphism given by 
\begin{align*}
M&=\begin{pmatrix} 0_n & I_n \\ -I_n & 0_n\end{pmatrix}.
\end{align*} 
\paragraph{Automorphisms of $\SU_{n}(q)$}
All linear outer automorphisms of $\SU_{n}(q)$ are induced by $\GU_n(q)$ the group of linear isometries of the hermitian form,  or are induced by $\Aut(\FF_{q^2})$ as for $\SL_n(q^2)$ with respect to an orthonormal basis.
The group $\Aut(\FF_{q^2})$ has order $2e$, where $q=p^e$, $p$ prime. 
We let $\GammaU_n(q)=\GU_n\rtimes \Aut(\FF_{q^2})$ and let $\PGammaU_n(q)$ denote its quotient over the center (consisting of the scalar matrices).
In this case, the transpose-inverse map $\trin$ with respect to a hyperbolic basis is the composition of the inner automorphism given by 
 \begin{align*}
M&=\begin{pmatrix} 0_n & I_n \\ I_n & 0_n\end{pmatrix}.
\end{align*} 
and the field automorphism $x\mapsto \bar{x}=x^q$ (with respect to the hyperbolic basis).

\paragraph{The group $\widehat{\Aut}(\FF_{q^2})$ of field automorphisms of $\SU_n(q)$ on a hyperbolic basis}
For $\GammaU_{2n}(q)$ note that $\Aut(\FF_{q^2})=\langle\alpha\rangle$ acts with respect to an orthonormal basis $\cU=\{u_1,\ldots,u_{2n}\}$ for the $\FF_{q^2}$-vector space $V$ with $\sigma$-hermitian form $\hform$ preserved by the group (see~\cite{Wil2009}).
.
We now identify a complement $\widehat{\Aut}(\FF_{q^2})$ of semilinear automorphisms of $\GU_{2n}(q)$ in $\GammaU_{2n}(q)$ with respect to a hyperbolic basis.
Fix the standard hyperbolic basis $\cH=\{e_i,f_i\colon i =1,2,\ldots,n\}$ so that the elements in $\GU(V,\hform)$ are represented by a matrix in $\GU_{2n}(q)$ with respect to $\cH$.
Let $\alpha\in \Aut(\FF_{q^2})$ act on $V$ via $\cU$. Then, $\cH^\alpha=\{e_i^\alpha,f_i^\alpha\colon i=1,2\ldots,n\}$ is also a hyperbolic basis for $V$, so for some $A\in \GU_{2n}(q)$, we have $A\cH=\cH^\alpha$. Now the composition $\hat{\alpha}=A^{-1}\after\alpha$ is an $\alpha$-semilinear map that fixes $\cH$. The corresponding automorphism of $\GU_{2n}(q)$ acts by applying $\alpha$ to the matrix entries.

\bre\label{rem:SU2 to SL2 Autos}
The following special case will be of particular interest when considering a Curtis-Tits standard pairs of type $\twA_3(q)$. 
In this case the action of $\hat{\alpha}$ as above on $\SU_4(q)$ translates via the standard identification maps (see~Definition~\ref{dfn:standard CT identification map}) to actions on $\SL_2(q)$ and $\SL_2(q^2)$ as follows.
The action on $\SL_2(q^2)$ is the natural entry-wise field automorphism action. The action on $\SL_2(q)$ will be a product of the natural entrywise action of $\hat{\alpha}$ and a diagonal automorphism $\diag(f,1)$, where $f\in \FF_q$ is such that $\hat{\alpha}(\eta)=f\eta$. Note that $N_{\FF_{q}/\FF_p}(f)=-1$, so in particular $\sigma=\hat{\alpha}^e$ translates to (left) conjugation by $\diag(-1,1)$ only. 

\ere
\bde\label{dfn:check autos of SU4}
Since the norm is surjective, there exists $\zeta\in \FF_{q^2}$  such that $N_{\FF_{q^2}/\FF_q}(\zeta)=f^{-1}$.
We then have that $\diag(\zeta,\zeta,\zeta^{-q},\zeta^{-q})\in \GU_4(q)$ acts trivially on $\SL_2(q^2)$ and acts as left conjugation by $\diag(f^{-1},1)$ on $\SL_2(q)$. It follows that the composition {$\widetilde{\alpha}$} of $\hat{\alpha}$ and this diagonal automorphism acts entrywise as $\alpha$ on both $\SL_2(q)$ and $\SL_2(q^2)$.
We now define  
\begin{align*}
{\widetilde{\Aut}}(\FF_{q^2})=\langle \widetilde{\alpha}\rangle\le \Aut(\SU_4(q)).
\end{align*}
\ede

\ble\label{lem:Aut(SL)}(See ~\cite{SchVan28,Wil2009}.)
\begin{enumerate}
\item As $\Sp_2(q)=\SL_2(q)\cong\SU_2(q)$, we have 
\begin{align*}
\Aut(\Sp_2(q))=\Aut(\SL_2(q))=\PGammaL_2(q)\cong\PGammaU_2(q)=\Aut(\SU_2(q)). 
\end{align*}
\item
In higher rank we have 
\begin{align*}
\Aut(\SL_n(q))&=\PGammaL_n(q)\rtimes \langle \trin\rangle\\
\Aut(\Sp_{2n}(q))&=\PGammaSp_{2n}(q)\\
\Aut(\SU_{n}(q))&=\PGammaU_{n}(q)
\end{align*}
\end{enumerate}
\ele
%\bre
%
%\ere

\subsubsection{Some normalizers and centralizers}
\bco\label{cor:Levi centralizer in SL3}
Let $G=\SL_3(q)$.
Let $\varphi\colon \SL_2(q)\to G$ given by 
$A\mapsto \begin{pmatrix} 1 & 0 \\ 0 & A\end{pmatrix}$ and let $L=\im \varphi$.
Then, 
\begin{align*}
C_{\Aut(G)}(L) & = \langle \diag(a,b,b)\colon a,b\in \FF_q^*\rangle\rtimes \langle \theta \rangle.
\end{align*}
where $\theta=\trin\after c_\nu\colon X^\theta\mapsto {}^t(\nu^{-1} X \nu)^{-1}$ and $\nu=\begin{pmatrix} 1 &  &  \\
  & 0 & -1\\
 & 1 & 0\\
\end{pmatrix}$.
\eco
\bpf
This follows easily from the fact that $\Aut(G)\cong\PGammaL_3(q)$. 
Let $\trin^i\after \alpha\after c_g$, where $c_g$ denotes conjugation by $g\in \GL_3(q)$ and $\alpha\in \Aut(\FF_q)$.
Using transvection matrices from $L$ over the fixed field $\FF_q^\alpha$ one sees that if $i=0$, then $g$ must be of the form $\diag(a,b,b)$, and if $i=1$, then it must be of the form 
 $\diag(a,b,b) \nu$, for some $a,b\in \FF_q^*$.
 Then, if $\alpha\ne\id$, picking transvections from $L$ with a few entries in $\FF_q-\FF_q^\alpha$ one verifies that $\alpha$ must be the identity.
\epf

\section{Classification of Curtis-Tits amalgams}\label{sec:classification of CT amalgams}

\subsection{Fundamental root groups in Curtis-Tits standard pairs}\label{subsub:fundamental root groups in CT standard pairs} 
\ble\label{lem:Levi1 with root group in Levi2}
Let $q$ be a power of the prime $p$.
Suppose that $(\amgrpG,\amgrpG_1,\amgrpG_2)$ is a Curtis-Tits standard pair of type $\liediag(q)$ as in Subsection~\ref{subsec:standard CT pairs}.
For $\{i,j\}=\{1,2\}$, let $\cS_j=\Syl_p(\amgrpG_j)$.
\begin{enumerate}
\item 
There exist two groups $\amgrpX_j^{i,\vep}\in \cS_j$ ($\vep=+,-$) such that for any $\amgrpX\in \cS_j$ we have
\begin{align*}
\langle \amgrpG_i,\amgrpX\rangle\le \amgrpP_i^\vep & \mbox{ if and only if }\amgrpX=\amgrpX_j^{i,\vep} 
\end{align*}
where $\amgrpP_i^+$ and $\amgrpP_i^-$ are the two parabolic subgroups of $\amgrpG$ containing $\amgrpG_i$.
If $\amgrpX\ne \amgrpX_j^{i,\vep}$, then 
\begin{align*}
\langle \amgrpG_i,\amgrpX\rangle 
=
\begin{cases} 
(\amgrpG_i\times\amgrpG_i^x)\rtimes \langle x\rangle  & \mbox{ if }\liediag(q)=C_2(2),\\
\amgrpG & \mbox{ else. }
\end{cases}
\end{align*}
where in the $C_2(2)$ case $\amgrpX=\langle x\rangle$.

\item 
We can select the signs $\vep$ so that 
$\amgrpX_i^{j,\vep}$ commutes with $\amgrpX_j^{i,-\vep}$, but not with $\amgrpX_j^{i,\vep}$ and, in fact  
$\langle \amgrpX_i^{j,\vep}, \amgrpX_j^{i,\vep}\rangle$ is contained in the unipotent radical $\amgrpU_{i,j}^\vep$ of a unique Borel subgroup of $\amgrpG_{i,j}$, namely $\amgrpB_{i,j}^\vep=\amgrpP_i^\vep\cap\amgrpP_j^\vep$.
\end{enumerate}
\ele
\bpf
We first prove part 1.~by considering all cases.
\paragraph{\bf ${\mathbf A_2(q)}$, ${\mathbf q\ge 3}$}
View $\amgrpG=\SL_3(q)=\SL(V)$ for some $\FF_q$-vector space with basis $\{e_1,e_2,e_3\}$.
By symmetry we may assume that $i=1$ and $j=2$.
Let $\amgrpG_1$ (resp. $\amgrpG_2$) stabilize $\langle e_1,e_2\rangle$  and fix $ e_3$ (resp. stabilize $\langle e_2,e_3\rangle$) and fix $e_1$).
A root group in $\amgrpG_2$ is of the form $\amgrpX_v=\Stab_{\amgrpG_2}(v)$ for some $v\in\langle e_2,e_3\rangle$.
We let $\amgrpX_2^+=\amgrpX_{e_2}$ and $\amgrpX_2^-=\amgrpX_{e_3}$.
It clear that for $\vep=+$ (resp. $\vep=-$) $\langle \amgrpG_1,\amgrpX_2^\vep\rangle=\amgrpP^\vep$  is contained in (but not equal to) the parabolic subgroup stabilizing $\langle e_1,e_2\rangle$  (resp. $\langle e_3\rangle$).
Now suppose that $\amgrpX\in \cS_2$ is different from $\amgrpX_2^\vep$ ($\vep=+,-$) and $\amgrpX=\amgrpX_{\lambda e_2+ e_3}$ for some $\lambda\in \FF_q^*$.
Consider the action of a torus element $d=\diag(\mu, \mu^{-1}, 1)\in \amgrpG_1$
by conjugation on $\amgrpG_2$. Then $\amgrpX^d=\amgrpX_{\mu\lambda e_2+e_3}$.
 Since $|\FF_q|\ge 3$, $\amgrpX^d\ne \amgrpX$ for some $d$ and so we have 
 \begin{align}
\langle \amgrpG_i,\amgrpX\rangle\ge \langle \amgrpG_i,\amgrpX,\amgrpX^d\rangle =
\langle\amgrpG_i,\amgrpG_j\rangle = \amgrpG.\label{eqn:G1XXd}
 \end{align}
\paragraph{\bf ${\mathbf A_2(2)}$}
In this case $\cS_2 =\{\amgrpX_2^+,\amgrpX=\langle r\rangle, \amgrpX_2^-\}$, where 
 $r$ is the Coxeter element fixing $e_1$ and interchanging $e_2$ and $e_3$.
It follows that  $\amgrpG_1^r$ is the stabilizer of the subspace decomposition $\langle e_2\rangle\oplus\langle e_1,e_3\rangle$ and hence $\langle \amgrpG_1,\amgrpX\rangle=\amgrpG$.

\paragraph{\bf ${\mathbf C_2(q)}$, ${\mathbf q\ge 3}$, $\mathbf \amgrpX$ short root}
We use  the notation of Subsection~\ref{subsec:standard CT pairs}.
First, let $i=2$, $j=1$, let $\amgrpG_2\cong \Sp_2(q)\cong\SL_2(q)$ be the stabilizer of $e_1$ and $e_3$ and let $\amgrpG_1\cong\SL_2(q)$ be the derived subgroup of the stabilizer of the isotropic $2$-spaces $\langle e_1,e_2\rangle$ and $\langle e_3,e_4\rangle$.
Root groups in $\amgrpG_1$ are of the form
 $\amgrpX_{u,v}=\Stab_{\amgrpG_1}(u)\cap \Stab_{\amgrpG_1}(v)$, where 
  $u\in \langle e_1,e_2\rangle$ and $v\in\langle e_3,e_4\rangle$ are orthogonal.
Let $\amgrpX_1^+=\amgrpX_{e_1,e_4}$ and $\amgrpX_1^-=\amgrpX_{e_2,e_3}$.
It is easy to verify that for $\vep=+$ (resp. $\vep=-$) $\langle \amgrpG_2,\amgrpX_1^\vep\rangle=\amgrpP^\vep$  is contained in the parabolic subgroup stabilizing $\langle e_1\rangle $ (resp. $\langle e_3\rangle$).
Now let $\amgrpX=\amgrpX_{e_1+\lambda e_2, e_3-\lambda^{-1} e_4}$ for some $\lambda\in\FF_q^*$.
Consider the action of a torus element $d=\diag(1, \mu^{-1}, 1, \mu)\in \amgrpG_2$
by conjugation on $\amgrpG_1$. 
Then $\amgrpX^d=\amgrpX_{\langle e_1+ \lambda\mu e_2\rangle,e_3-\lambda^{-1}\mu^{-1} e_4}$.
 Since $q\ge 3$, $\amgrpX^d\ne \amgrpX$ for some $d$ and so, for $i=1$, and these $\amgrpG_2$, $\amgrpX$ and $d$, we have~\eqref{eqn:G1XXd} again.

\paragraph{\bf ${\mathbf C_2(q)}$, ${\mathbf q\ge 4}$, $\mathbf \amgrpX$ long root}
Now, we let $i=1$ and $j=2$.
Root groups in $\amgrpG_2$ are of the form $\amgrpX_u=\Stab_{\amgrpG_2}(u)$ where $u\in \langle e_2,e_4\rangle$.
Let $\amgrpX_2^+=\amgrpX_{e_2}$  and $\amgrpX_2^-=\amgrpX_{e_4}$.
It is easy to verify that for $\vep=+$ (resp. $\vep=-$) $\langle \amgrpG_1,\amgrpX_2^\vep\rangle=\amgrpP^\vep$  is contained in the parabolic subgroup stabilizing $\langle e_1,e_2\rangle $ 
(resp. $\langle e_1,e_4\rangle$).
Now let $\amgrpX=\amgrpX_{e_2+\lambda e_4}$ for some $\lambda\in\FF_q^*$.

Consider the action of a torus element $d=\diag(\mu, \mu^{-1}, \mu^{-1}, \mu)\in \amgrpG_1$
by conjugation on $\amgrpG_2$. 
Then $\amgrpX^d=\amgrpX_{\mu e_2+ \mu^{-1}\lambda e_4}$.
Now if $q\ge 4$, then $\amgrpX^d\ne \amgrpX$ for some $d$ and so  so, for these $\amgrpG_1$, $\amgrpX$ and $d$, we have~\eqref{eqn:G1XXd} again.

\paragraph{\bf ${\mathbf C_2(q)}$, ${\mathbf q =3}$, $\mathbf \amgrpX$ long root}
The proof for the case $q\ge 4$ does not yield the result since, for $q=3$, the element $d$ centralizes $\amgrpG_2$.
A direct computation in GAP shows that the conclusion still holds, though.
Let $x\in \amgrpX=\amgrpX_{e_2+e_4}$ send $e_2$ to $e_4$.
Then $\amgrpG_1$ and $\amgrpG_1^x$ contains two short root groups fixing $e_1$ and $e_3$.
Their commutators generate a long root group fixing $e_1$, $e_2$, and $e_4$, while being transitive on the points $\langle e_3+\lambda e_1\rangle$. Further conjugation with an element in $\amgrpG_1$ interchanging the points $\langle e_1\rangle$ and $\langle e_2\rangle$ yields a long root group in $\amgrpG_2$ different from $\amgrpX$ and we obtain an equation like~\eqref{eqn:G1XXd} again.

\paragraph{\bf ${\mathbf C_2(2)}$} 
First note that $\amgrpG\cong\Sp_4(2)\cong O_5(2)$ is self point-line dual, so we only need to consider the case where $\amgrpG_2=\Stab_\amgrpG(e_1)\cap\Stab_\amgrpG(e_3)$ and 
 $\amgrpG_1=\Stab_\amgrpG(\langle e_1,e_2\rangle)\cap \Stab_\amgrpG(\langle e_3,e_4\rangle)$.
 Now $\cS_1=\{\amgrpX_1^+, \amgrpX_1^-, \langle x\rangle\}$, where $x$ is the permutation matrix of $(1,2)(3,4)$.
The conclusion follows easily.

\paragraph{\bf ${\mathbf\twA_3(q)}$}
We use  the notation of Subsection~\ref{subsec:standard CT pairs}.
First, let $i=2$, $j=1$, let $\amgrpG_2\cong \SU_2(q)\cong\SL_2(q)$ be the stabilizer of $e_1$ and $e_3$ and let $\amgrpG_1\cong\SL_2(q^2)$ be the derived subgroup of the simultaneous stabilizer in $\amgrpG=\SU_4(q)$ of the isotropic $2$-spaces $\langle e_1,e_2\rangle$ and $\langle e_3,e_4\rangle$.
Root groups in $\amgrpG_1$ are of the form
 $\amgrpX_{u,v}=\Stab_{\amgrpG_1}(u)\cap \Stab_{\amgrpG_1}(v)$, where 
  $u\in \langle e_1,e_2\rangle$ and $v\in\langle e_3,e_4\rangle$ are orthogonal.
Let $\amgrpX_1^+=\amgrpX_{e_1,e_4}$ and $\amgrpX_1^-=\amgrpX_{e_2,e_3}$.
It is easy to verify that for $\vep=+$ (resp. $\vep=-$) $\langle \amgrpG_2,\amgrpX_1^\vep\rangle=\amgrpP^\vep$  is contained in the parabolic subgroup stabilizing $\langle e_1\rangle $ (resp. $\langle e_3\rangle$).
Now let $\amgrpX=\amgrpX_{e_1+\lambda e_2, e_3-\lambda^{-\sigma} e_4}$ for some $\lambda\in\FF_{q^2}^*$.
Consider the action of a torus element $d=\diag(1, \mu^{-1}, 1, \mu)\in \amgrpG_2$ 
(with $\mu\in \FF_q^*$)
by conjugation on $\amgrpG_1$. 
Then $\amgrpX^d=\amgrpX_{e_1+ \lambda\mu e_2 ,e_3-\lambda^{-\sigma}\mu^{-1} e_4}$.
There are $q-1$ choices for $\mu$, so if $q\ge 3$, then $\amgrpX^d\ne \amgrpX$ for some $d$. Hence, for $i=1$, and these $\amgrpG_2$, $\amgrpX$ and $d$, we have~\eqref{eqn:G1XXd} again.
The case $q=2$ is a quick GAP calculation.

Now, we let $i=1$ and $j=2$.
Root groups in $\amgrpG_2$ are of the form $\amgrpX_u=\Stab_{\amgrpG_2}(u)$ where $u\in \langle e_2,e_4\rangle$ is isotropic.
Let $\amgrpX_2^+=\amgrpX_{e_2}$  and $\amgrpX_2^-=\amgrpX_{e_4}$.
It is easy to verify that for $\vep=+$ (resp. $\vep=-$) $\langle \amgrpG_1,\amgrpX_2^\vep\rangle=\amgrpP^\vep$  is contained in the parabolic subgroup stabilizing $\langle e_1,e_2\rangle $ 
(resp. $\langle e_1,e_4\rangle$).
Now let $\amgrpX=\amgrpX_{e_2+\lambda e_4}$ for some $\lambda\in\FF_{q^2}^*$ where 
 $\Tr(\lambda)=\lambda+\lambda^\sigma=0$.
Consider the action of a torus element $d=\diag(\mu, \mu^{-1}, \mu^{-\sigma}, \mu^\sigma)\in \amgrpG_1$ (for some $\mu\in \FF_{q^2}^*$) by conjugation on $\amgrpG_2$. 
Then $\amgrpX^d=\amgrpX_{\mu e_2+ \mu^{-\sigma}\lambda e_4}$.
The $q^2-1$ choices for $\mu$ result in $q-1$ different conjugates. 
Thus, if $q-1\ge 2$, then $\amgrpX^d\ne \amgrpX$ for some $d$ and so  so, for these $\amgrpG_1$, $\amgrpX$ and $d$, we have~\eqref{eqn:G1XXd} again.
The case $q=2$ is a quick GAP calculation.
Namely, in this case, $\amgrpX=\langle x\rangle$, where $x$ is the only element of order $2$ in $\amgrpG_2\cong S_3$ that does not belong to $\amgrpX_1^+\cup \amgrpX_1^-$; it is the Coxeter element that fixes $e_1$ and $e_3$ and interchanges $e_2$ and $e_4$.
Now $\langle \amgrpG_1,\amgrpG_1^r\rangle$ contains the long root group generated by the commutators of the short root group fixing $e_1$ in $\amgrpG_1$ and $\amgrpG_1^r$, and likewise for $e_2$, $e_3$, and $e_4$.
In particular, we have 
\begin{align}
\langle \amgrpG_1, \amgrpX\rangle\ge \langle \amgrpG_1,\amgrpG_1^r\rangle\ge \langle\amgrpG_1,\amgrpG_2\rangle=\amgrpG \label{eqn:<Gj,Gj^X> =G}
\end{align}

We now address part 2.
Note that the positive and negative fundamental root groups with respect to the torus $\amgrpB_{i,j}^+\cap \amgrpB_{i,j}^-$ satisfy the properties of $\amgrpX_i^{j,\vep}$ and  $\amgrpX_j^{i,\vep}$ so by the uniqueness statement in 1.~they must be equal.
Now the claims in part 2.~are the consequences of the Chevalley commutator relations.
\epf

\bre
 \ere

Explicitly, the groups $\{\amgrpX_i^+,\amgrpX_i^-\}$ ($i=1,2$), possibly up to a switch of signs, for the Curtis-Tits standard pairs are as follows.

For $\liediag=A_2$, we have
 \begin{align*}
\amgrpX_1^+&=\left\{\begin{pmatrix} 
1 & b  &  \\
 0 & 1  &  \\
 &  & 1  \end{pmatrix}\colon b\in\FF_{q} \right\},
  \mbox{ and  }
\amgrpX_1^-=\left\{\begin{pmatrix} 
1 &  0 & \\
 c & 1   & \\
 &  & 1  \end{pmatrix}\colon c\in\FF_{q} \right\}, 
  \\
\amgrpX_2^+&=\left\{\begin{pmatrix} 
1 &   &  \\
  & 1  &  b\\
 &  0 & 1  \end{pmatrix}\colon b\in\FF_{q} \right\},
  \mbox{ and  }
\amgrpX_2^-=\left\{\begin{pmatrix} 
1 &   & \\
  & 1   & 0\\
 &  c & 1  \end{pmatrix}\colon c\in\FF_{q} \right\}.
  \end{align*}

For $\liediag=C_2$, we have
 \begin{align*}
\amgrpX_1^+&=\left\{\begin{pmatrix} 
1 & b & & \\
0 & 1  & & \\
 &  & 1 & 0 \\
 &  &  -b & 1 \end{pmatrix}\colon b \in\FF_{q} \right\},
\mbox{ and }
\amgrpX_1^-=\left\{\begin{pmatrix} 
1 & 0 & & \\
c & 1  & & \\
 &  & 1 & -c \\
 &  &  0 & 1 \end{pmatrix}\colon c \in\FF_{q} \right\},\\
\amgrpX_2^+&=\left\{\begin{pmatrix} 
1 &  & & \\
 & 1  & & b\\
 &  & 1& \\
 & 0 & &1 \end{pmatrix}\colon b\in\FF_{q} \right\},
  \mbox{ and  }
\amgrpX_2^-=\left\{\begin{pmatrix} 
1 &  & & \\
 & 1  & & 0\\
 &  & 1& \\
 & c & &1 \end{pmatrix}\colon c\in\FF_{q} \right\}.
  \end{align*}

For $\liediag=\twA_3$, we have (with $\eta\in \FF_{q^2}$ of trace $0$), 
 \begin{align*}
\amgrpX_1^+&=\left\{\begin{pmatrix} 
1 & b & & \\
0 & 1  & & \\
 &  & 1 & 0 \\
 &  &  -b & 1 \end{pmatrix}\colon b \in\FF_{q^2} \right\},
\mbox{ and }
\amgrpX_1^-=\left\{\begin{pmatrix} 
1 & 0 & & \\
c & 1  & & \\
 &  & 1 & -c \\
 &  &  0 & 1 \end{pmatrix}\colon c \in\FF_{q^2} \right\},\\
 \amgrpX_2^+&=\left\{\begin{pmatrix} 
1 &  & & \\
 & 1  & & b\eta\\
 &  & 1& \\
 & 0 & &1 \end{pmatrix}\colon b\in\FF_{q} \right\},
  \mbox{ and  }
\amgrpX_2^-=\left\{\begin{pmatrix} 
1 &  & & \\
 & 1  & & 0\\
 &  & 1& \\
 & c\eta^{-1} & &1 \end{pmatrix}\colon c\in\FF_{q} \right\}. 
 \end{align*} 
\subsection{Weak systems of fundamental groups}\label{subsub:weak systems} 
In this subsection we show that a Curtis-Tits amalgam with $3$-spherical diagram determines 
 a collection of subgroups of the vertex groups, called a weak system of fundamental root groups.
We then use this to determine the coefficient system of the amalgam in the sense of~\cite{BloHof2013}, which, in turn is applied to classify these amalgams up to isomorphism.

\bde\label{dfn:weak system of fundamental root groups}
Suppose that $\amG=\{\amgrpG_i,\amgrpG_{i,j},\amg_{i,j}\mid i,j\in I\}$ is a  CT amalgam.
For each $i\in I$ let $\amgrpX_i^+, \amgrpX_i^-\le \amgrpG_i$ be a pair of opposite root groups.
We say that $\{\amgrpX_i^+, \amgrpX_i^-\mid i\in I\}$ is a {\dfn weak system of fundamental root groups} if, for any edge $\{i,j\}\in \edg$ there are opposite Borel groups $\amgrpB_{i,j}^+$  and $\amgrpB_{i j}^-$ in $\amgrpG_{i,j}$, each of which contains exactly one of $\{\bamgrpX_i^+,\bamgrpX_i^-\}$.

We call $\amG$ {\em orientable} if we can select  $\amgrpX_i^\vep$, $\amgrpB_{i j}^\vep$ ($\vep=+,-$) for all $i,j\in \vrtc$ such that 
 $\bamgrpX_i^\vep, \bamgrpX_j^\vep\le \amgrpB_{ i j}^\vep$.
If this is not possible, we call $\amG$ {\em non-orientable}.
 \ede

The relation between root groups and Borel groups is given by the following well-known fact.
\ble\label{lem:B = N_G(X)}
Let $q$ be a power of the prime $p$.
Let $\amgrpG$ be a universal group of Lie type $\liediag(q)$
and let $\amgrpX$ be a Sylow $p$-subgroup.
Then, $N_\amgrpG(\amgrpX)$ is the unique Borel group $\amgrpB$ of $\amgrpG$ containing $\amgrpX$.
\ele 

\bpr\label{prop:existence of weak system of fundamental root groups}
Suppose that $\amG=\{\amgrpG_i,\amgrpG_{i,j},\amg_{i,j}\mid i,j\in I\}$ is a  CT amalgam with connected $3$-spherical diagram $\liediag$.
If $\amG$ has a non-trivial completion $(\compG,\compg)$, then it has a unique weak system of fundamental root groups.
\epr
\bpf
We first show that there is some weak system of fundamental root groups.
For every edge $\{i,j\}$, let $\amgrpX_j^{i,\vep}$ be the groups of Lemma~\ref{lem:Levi1 with root group in Levi2}.
Suppose that there is some subdiagram $\liediag_J$ with $J=\{i,j,k\}$ in which $j$ is connected to both $i$ and $k$, such that $\{\amgrpX_j^{i,+},\amgrpX_j^{i,-}\}\ne\{\amgrpX_j^{k,+},\amgrpX_j^{k,-}\}$ as sets.
Without loss of generality assume that $\liediag_{i,j}=A_2$ (by $3$-sphericity) and moreover, that $\amgrpX_j^{k,+}\not\in \{\amgrpX_j^{i,+},\amgrpX_j^{i,-}\}$.
For any subgroup $\amgrpH$ of a group in $\amG$, write $\compH=\compg(\amgrpH)$.
Now note that $\compX_k^{j,-}$ commutes with $\compX_j^{k,+}$  and since $\liediag$ contains no triangles it also commutes with $\compG_i$. But then $\compX_k^{j,-}$ commutes with $\langle \compX_j^{k,+},\compG_i\rangle$ which, by Lemma~\ref{lem:Levi1 with root group in Levi2}, equals $\compG_{i,j}$
 (this is where we use that $\liediag_{i,j}=A_2$), contradicting that 
 $\compX_k^{j,-}$ does not commute with $\compX_j^{k,-}\le \compG_{i,j}$.
Thus, if there is a completion, then by connectedness of $\liediag$, for each $i\in I$ we can pick a $j\in I$ so that $\{i,j\}\in E$ 
 and set $\amgrpX_i^\pm=\amgrpX_i^{j,\pm}$ and drop the superscript.
We claim that $\{\amgrpX_i^\pm\mid i\in I\}$ is a weak system of fundamental root groups.
But this follows from part 2.~of Lemma~\ref{lem:Levi1 with root group in Levi2}.

The uniqueness derives immediately from the fact that by Lemma~\ref{lem:Levi1 with root group in Levi2}, $\amg_{i,j}(\amgrpX_j^{i,+})$ and $\amg_{j,i}(\amgrpX_j^{i,-})$ are the only two Sylow $p$-subgroups in $\amg_{j,i}(\amgrpG_j)$ which do not generate $\amgrpG_{i,j}$ with $\amg_{i,j}(\amgrpG_i)$.
\epf

\medskip
An immediate consequence of the results above is the following observation.
\bco\label{cor:edge autos fix or swap + and -}
Suppose that $\amG=\{\amgrpG_i,\amgrpG_{i,j},\amg_{i,j}\mid i,j\in I\}$ is a  CT amalgam with connected $3$-spherical diagram $\liediag$.
Then, 
an element of $N_{\Aut(\amgrpG_{i,j})}(\bamgrpG_i,\bamgrpG_j)$ either fixes
 each of the pairs $(\bamgrpX_i^+,\bamgrpX_i^-)$, $(\bamgrpX_j^+,\bamgrpX_j^-)$, and  $(\amgrpB_{i,j}^+,\amgrpB_{i,j}^-)$ or it reverses each of them.
In particular,  
\begin{align*}
N_{\Aut(\amgrpG_{i,j})}(\bamgrpG_i,\bamgrpG_j) &  =  N_{\Aut(\amgrpG_{i,j})}(\{\bamgrpX_i^+,\bamgrpX_i^-\})\cap N_{\Aut(\amgrpG_{i,j})}(\{\bamgrpX_j^+,\bamgrpX_j^-\}).
\end{align*} 
\eco

\subsection{The coefficient system of a  Curtis-Tits amalgam}\label{subsec:coefficient system}
The automorphisms of a Curtis-Tits standard pair will be crucial in the classification of  Curtis-Tits amalgams and we will need some detailed description of them.

We now fix a Curtis-Tits amalgam  $\famG=\{\amgrpG_i,\amgrpG_{i,j},\famg_{i,j}\mid i,j\in I\}$ of type $\liediag(q)$, where for every $i,j\in I$, $\famg_{i,j}$ is the standard identification map of Definition~\ref{dfn:standard CT identification map}.
Then,  $\famG$ has a weak system of fundamental root groups $\Chi=\{\{\amgrpX_i^+,\amgrpX_i^-\}\colon i\in I\}$ as in Subsection~\ref{subsub:fundamental root groups in CT standard pairs}.

\bre\label{rem:all  Curtis-Tits amalgams have same type}
Let  $\amG=\{\amgrpG_i,\amgrpG_{i,j},\amg_{i,j}\mid i,j\in I\}$ be a Curtis-Tits amalgam over $\FF_q$ with given diagram $\liediag$. 
Next suppose that $\liediag$ is connected $3$-spherical, and that $\famG$ and $\amG$ are non-collapsing.
Then, by Proposition~\ref{prop:existence of weak system of fundamental root groups}, $\famG$ and $\amG$ each have a weak system of fundamental root groups. Now note that for each $i\in I$, $\Aut(\amgrpG_i)$ is $2$-transitive on the set of Sylow $p$-subgroups.
Thus, for each $i\in I$ and all $j\in I-\{i\}$, we can replace $\amg_{i,j}$ by $\amg_{i,j}\after\alpha_i$, to form a new amalgam isomorphic to $\amG$, whose weak system of fundamental root groups is exactly $\Chi$.
Thus, in order to classify non-collapsing Curtis-Tits amalgams over $\FF_q$ with diagram $\liediag$ up to isomorphism, it suffices to classify those whose weak system of fundamental root groups is exactly $\Chi$.
\ere

\bde\label{dfn:coefficient system}
Suppose that $\amG=\{\amgrpG_i,\amgrpG_{i,j},\amg_{i,j}\mid i,j\in I\}$ is a  Curtis-Tits amalgam over $\FF_q$ with connected $3$-spherical diagram $\liediag$.
Denote the associated weak system of fundamental root groups as $\Chi=\{\{\amgrpX_i^+,\amgrpX_i^-\}\colon i\in I\}$.
The {\dfn coefficient system associated to $\amG$} is the collection 
 $\amA=\{\amgrpA_i,\amgrpA_{i,j},\ama_{i,j}\mid i,j\in I\}$ where, for any $i,j\in I$ we set 
 \begin{align*}
 \amgrpA_i&=N_{\Aut(\amgrpG_i)}(\{\amgrpX_i^+,\amgrpX_i^-\}), \\
 \amgrpA_{i,j}&=N_{\Aut(\amgrpG_{i,j})}(\{\bamgrpX_i^\vep \colon\vep=+,-\})\cap N_{\Aut(\amgrpG_{i,j})}(\{\bamgrpX_j^\vep \colon\vep=+,-\}), \\
 \ama_{i,j}&\colon \amgrpA_{i,j}\to\amgrpA_j \mbox{ is given by restriction: } \varphi\mapsto \amg_{j,i}^{-1}\after \rho_{i,j}(\varphi)\after\amg_{j,i}.
  \end{align*}
where $\rho_{i,j}(\varphi)$ is the restriction of $\varphi$ to $\bamgrpG_j\le \amgrpG_{i,j}$.
\ede

From now on we let $\amA$ be the coefficient system associated to $\famG$.
The significance for the classification of  Curtis-Tits amalgams with weak system of fundamental root groups is as follows:
\bpr\label{prop:coefficient systems}
Suppose that $\amG$ and $\amG^+$ are  Curtis-Tits amalgams with diagram $\liediag$ over $\FF_q$ with weak system of fundamental root groups $\Chi$.
\begin{enumerate}
\item For all $i,j\in I$, we have $\amg_{i,j}=\famg_{i,j}\after\delta_{i,j}$ and $\amg_{i,j}^+=\famg_{i,j}\after\delta_{i,j}^+$ for some $\delta_{i,j}, \delta_{i,j}^+\in \amgrpA_i$,
\item For any isomorphism $\phi\colon \amG\to\amG^+$ and $i,j\in I$, we have $\phi_i\in \amgrpA_i$, $\phi_{\{i,j\}}\in \amgrpA_{i,j}$, and $\ama_{i,j}(\phi_{\{i,j\}})=\delta_{i,j}^+\after\phi_i\after\delta_{i,j}^{-1}$.
\end{enumerate}
\epr
\bpf
Part 1.~follows since, for any $i,j\in I$ we have $\amg_{i,j}^{-1}\after\famg_{i,j}\in \Aut(\amgrpG_1)$ and  
\begin{align*}
\{\amg_{i,j}(\amgrpX_i^+),\amg_{i,j}(\amgrpX_i^-)\}=\{\famg_{i,j}(\amgrpX_i^+),\famg_{i,j}(\amgrpX_i^-)\}.
\end{align*}
Part 2.~follows from Corollary~\ref{cor:edge autos fix or swap + and -} since, for any $i,j\in I$,  
\begin{align*}
(\amgrpG_{i,j},\amg_{i,j}(\amgrpG_i),\amg_{j,i}(\amgrpG_j))
=(\amgrpG_{i,j},\famg_{i,j}(\amgrpG_i),\famg_{j,i}(\amgrpG_j)) =(\amgrpG_{i,j},\amg^+_{i,j}(\amgrpG_i),\amg^+_{j,i}(\amgrpG_j)).
\end{align*}
\epf

We now determine the groups appearing in the coefficient system $\amA$ associated to $\famG$.

\ble\label{lem:structure of CT A_i groups}
Fix $i\in I$ and let $q$ be such that $\amgrpG_i\cong\SL_2(q)$.
Then, 
\begin{align*}
\amgrpA_i&=\amgrpT_i\rtimes \amgrpC_i,
\end{align*}
where $\amgrpT_i$ is the subgroup of diagonal automorphisms in $\PGL_2(q)$ and $\amgrpC_i=\langle \trin, \Aut(\FF_q)\rangle$.
\ele
\bpf
This follows from the fact that via the standard embedding map $\famg_{i,j}$ the groups $\amgrpX_i^+$  and $\amgrpX_i^-$ of the weak system of fundamental root groups are the subgroups of unipotent upper and lower triangular matrices in $\SL_2(q)$.
\epf
\ble\label{lem:coefficient system groups}
Let $\amA$ be the coefficient system associated to the standard Curtis-Tits amalgam $\famG$ of type $\liediag(q)$ and the weak system of fundamental root groups $\Chi$.

If $\Gamma=A_1\times A_1$, we have $\amgrpG_{i,j}=\amgrpG_i\times\amgrpG_j$, $\famg_{i,j}$ and $\famg_{j,i}$ are identity maps, and  
\begin{align}
\amgrpA_{i,j}&=\amgrpA_i\times\amgrpA_j\cong \amgrpT_{i,j}\rtimes\amgrpC_{i,j}.\label{eqn:CT N Xi Xj A1timesA1}
\end{align}
where $\amgrpT_{i,j}=\amgrpT_i\times\amgrpT_j$ and $\amgrpC_{i,j}=\amgrpC_i\times\amgrpC_j$.
Otherwise, 
\begin{align*}
\amgrpA_{i,j}&=\amgrpT_{i,j}\rtimes\amgrpC_{i,j},
\end{align*}
where 
\begin{align}
\amgrpC_{i,j}&=\begin{cases}
\Aut(\FF_q)\times\langle \trin\rangle  & \mbox{ for }\liediag=A_2, C_2\\
\widetilde{\Aut}(\FF_{q^2})\times\langle \trin\rangle & \mbox{ for }\liediag=\twA_3 
\end{cases}
\end{align}
and  $\amgrpT_{i,j}$ denotes the image of the standard torus $\GD$ in $\Aut(\amgrpG_{i,j})$.
Note that 
\begin{align*}
\GD&=\begin{cases}
\langle \diag(a,b,c)\colon a,b,c\in \FF_q^*\rangle\le\GL_3(q) & \mbox{ if } \liediag=A_2\\
\langle \diag(ab,a^{-1}b,a^{-1},a)\colon a,b\in\FF_q^*\rangle\le\GSp_4(q) & \mbox{ if } \liediag=C_2\\ 
\langle \diag(a,b,a^{-q},b^{-q})\colon a,b\in \FF_{q^2}^* \rangle\le\GU_4(q) &
\mbox{ if } \liediag=\twA_3.\\ 
\end{cases}
\end{align*}
\ele
\bre
Remarks on Lemma~\ref{lem:coefficient system groups}
\begin{enumerate}
\item We view $\Sp_{2n}(q)$ and $\SU_{2n}(q)$ as a matrix group with respect to a symplectic (resp.~hyperbolic) basis for the $2n$-dimensional vector space $V$ and $\Aut(\FF_q)$ (resp. $\widehat{\Aut}(\FF_{q^2})$) acts entrywise on the matrices.
\item The map $\trin$ is the transpose-inverse map of Subsection~\ref{subsec:automorphism groups}.
\item Recall that in the $\twA_3$ case, Remark~\ref{rem:SU2 to SL2 Autos}~and~Definition~\ref{dfn:check autos of SU4} describe the actions of $\widetilde{\Aut}(\FF_{q^2})\le \amgrpC_{i,j}$ on $\amgrpG_i$ and $\amgrpG_j$ via the standard identification maps.
\end{enumerate}
\ere

\bpf
We first consider the $A_1\times A_1$ case of~\eqref{eqn:CT N Xi Xj A1timesA1}.
When $\liediag=A_1\times A_1$, then $\amgrpG_{i,j}=\bamgrpG_i\times\bamgrpG_j$ and since the standard root groups 
 $\bamgrpX_i^\pm$ generate $\bamgrpG_i$ ($i=1,2$), their simultaneous normalizer must also normalize $\bamgrpG_i$ and $\bamgrpG_j$. Thus the claim follows from Lemma~\ref{lem:structure of CT A_i groups}.

We now deal with all remaining cases simultaneously.
In the $\twA_3$ case we note that from Remark~\ref{rem:SU2 to SL2 Autos} and Definition~\ref{dfn:check autos of SU4} we see that $\widetilde{\Aut}(\FF_{q^2})\le \amgrpT_{i,j}\rtimes \widehat{\Aut}(\FF_{q^2})$
 is simply a different complement to $\amgrpT_{i,j}$, so it suffices to prove the claim with $\widetilde{\Aut}(\FF_{q^2})$ replaced by $\widehat{\Aut}(\FF_{q^2})$.

Consider the descriptions of the set $\{\amgrpX_i^+,\amgrpX_i^-\}$ in all cases from Subsection~\ref{subsub:fundamental root groups in CT standard pairs}.
We see that since $\trin$ acts by transpose-inverse, it interchanges $\amgrpX_i^+$ and $\amgrpX_i^-$ for $i=1,2$ in all cases, hence it also interchanges positive and negative Borel groups (see Corollary~\ref{cor:edge autos fix or swap + and -}).
Thus it suffices to consider those automorphisms that normalize the positive and negative fundamental root groups.
Since all field automorphisms (of $\Aut(\FF_q)$ and $\widehat{\Aut}(\FF_{q^2})$) act entrywise, they do so. Clearly so does $\GD$.
Thus we have established $\supseteq$.

We now turn to the reverse inclusion. By Lemma~\ref{lem:Aut(SL)} and the description of the automorphism groups in Subsection~\ref{subsec:automorphism groups} any automorphism of $\amgrpG_{i,j}$ is a product of the form $g\alpha\trin^i$ where $g$ is linear, $\alpha$ is a field automorphism (from $\widehat{\Aut}(\FF_{q^2})$ in the $\twA_3$ case)  and $i=0,1$.
As we saw above $\trin$ and $\alpha$ preserve the root groups, so it suffices to describe $g$ in case it preserves the sets of opposite root groups.
A direct computation shows that $g$ must be in $\GD$.
\epf

\medskip

Next we describe the connecting maps $\ama_{i,j}$ of $\amA$.

\ble\label{lem: coefficient system connecting maps}
Let $\amA$ be the coefficient system of the standard Curtis-Tits amalgam $\famG$ over $\FF_q$ with diagram $\liediag$ and weak system of fundamental root groups $\Chi$.
Fix $i,j\in I$ and let $(\amgrpG_{i,j},\bamgrpG_i,\bamgrpG_j)$ be a Curtis-Tits standard pair in $\famG$ with diagram $\liediag_{i,j}$. Denote $\ama=(\ama_{j,i},\ama_{i,j})\colon \amgrpA_{i,j}\to \amgrpA_i\times\amgrpA_j$.
Then, we have the following:
\begin{enumerate}
\item If $\liediag_{i,j}=A_1\times A_1$, then $\ama$ is an isomorphism inducing $\amgrpT_{i,j}\cong\amgrpT_i\times\amgrpT_j$ and $\amgrpC_{i,j}\cong\amgrpC_i\times\amgrpC_j$.
\item If $\liediag_{i,j}=A_2$, or $\twA_3$, then $\ama\colon \amgrpT_{i,j}\to \amgrpT_i\times\amgrpT_j$ is bijective.
\item If $\liediag_{i,j}=C_2$, then $\ama\colon \amgrpT_{i,j}\stackrel{\cong}{\to}\amgrpT_i^2\times\amgrpT_j$ has index $1$ or $2$ in $\amgrpT_i\times\amgrpT_j$ depending on whether $q$ is even or odd.
\item If $\liediag_{i,j}=A_2$ or $C_2$, then 
 $\ama\colon\amgrpC_{i,j}\to\amgrpC_i\times\amgrpC_j$ is given by  $\trin^s\alpha\mapsto (\trin^s\alpha,\trin^s\alpha)$ (for $s\in \{0,1\}$ and $\alpha\in \Aut(\FF_q)$) which is a diagonal  embedding.
\item If $\liediag=\twA_3$, then $\ama\colon {\amgrpC}_{i,j}\to \amgrpC_i\times\amgrpC_j$, is given by 
 $\trin^s\widetilde{\alpha}^r\mapsto (\trin^s\alpha^r,\trin^s\alpha^r)$ (for $s\in \{0,1\}$, $r\in \NN$, and 
  $\alpha\colon x\mapsto x^p$ for $x\in \FF_{q^2}$.  Here $\widetilde{\sigma}\mapsto (\sigma,\id)$.
\end{enumerate}
\ele
\bre
\begin{enumerate}
\item In 4.~$\trin$ acts as transpose-inverse and $\alpha$ acts entry-wise on $\amgrpG_{i,j}$, $\amgrpG_i$ and $\amgrpG_j$.
\end{enumerate}
\ere
\bpf
1.~ This is immediate from Lemma~\ref{lem:coefficient system groups}.

For the remaining cases, recall that for any $\varphi\in \amgrpA_{i,j}$, we have 
 $\ama_{i,j}\colon \varphi\mapsto \famg_{j,i}^{-1}\after \rho_{i,j}(\varphi)\after\famg_{j,i}$, 
where $\rho_{i,j}(\varphi)$ is the restriction of $\varphi$ to $\bamgrpG_j\le \amgrpG_{i,j}$ (Definition~\ref{dfn:coefficient system}) and $\famg_{i,j}$ is the standard identification map of Definition~\ref{dfn:standard CT identification map}.
Note that for $\liediag_{i,j}=A_2,C_2$ the standard identification map transforms the automorphism $\rho_{j,i}(\varphi)$ of $\bamgrpG_i$ essentially to the ``same" automorphism $\varphi$ of $\amgrpG_i$, whereas for $\liediag_{i,j}=\twA_3$, we must take Remark~\ref{rem:SU2 to SL2 Autos} into account.

2.~
Let $\liediag_{i,j}=A_2$.
Every element of $\amgrpT_{i,j}$  ($\amgrpT_i$, and $\amgrpT_j$ respectively) is given by a unique matrix of the form $\diag(a,1,c)$ ($\diag(a,1)$, and $\diag(1,c)$), and  we have 
\begin{align*}
(\ama_{j,i},\ama_{i,j})\colon \diag(a,1,c)\mapsto (\diag(a,1),\diag(1,c)) && (a,c\in \FF_q^*),
\end{align*}
which is clearly bijective. 
In the $\twA_3$ case, every element of  $\amgrpT_{i,j}$  ($\amgrpT_i$, and $\amgrpT_j$ respectively) is given by a unique matrix of the form $\diag(ab^{-1},1,a^{-q}b^{-1},b^{-(q+1)})$ ($\diag(1,b^{-(q+1)})$, and $\diag(ab^{-1},1)$), and we have 
\begin{align*}
\ama\colon \diag(ab^{-1},1,a^{-q}b^{-1},b^{-(q+1)}){\mapsto} (\diag(1,b^{-(q+1)}),\diag(ab^{-1},1)).
\end{align*}
This map is onto since $N_{\FF_{q^2}/\FF_q}\colon b\mapsto b^{q+1}$ is onto. Its kernel is trivial, as it is given by pairs $(a,b)\in \FF_{q^2}$, with $a=b$ and $b^{q+1}=1$ so that also $a^{-q}b^{-1}=1$.

3.~In the $C_2$-case,  every element of $\amgrpT_{i,j}$ is given by a unique diagonal matrix $\diag(a^2b, b,1,a^2)$ ($a,b\in \FF_q$).
Every element of $\amgrpT_i$ (resp. $\amgrpT_j$) is given by a unique $\diag(c,1)$ (resp. $\diag(d,1)$).
Now we have
\begin{align*}
\ama&\colon\diag(a^2b, b,1,a^2)\mapsto (\diag(a^2,1), \diag(ba^{-2},1)).
\end{align*}
It follows that $\ama$ is injective and has image 
 $\amgrpT_i^2\times\amgrpT_j$.
The rest of the claim follows.

4.~The field automorphism $\alpha\in \Aut(\FF_q)$ acts entrywise on the matrices in $\amgrpG_{i,j}=\SL_3(q)$, or  $\Sp_4(q)$, and $\amgrpG_i=\amgrpG_j=\SL_2(q)$.   In the case $\amgrpG=\Sp_4(q)$, we saw in Subsection~\ref{subsec:automorphism groups} that $\trin$ is inner and coincides with conjugation by $M$. This clearly restricts to conjugation by $\mu$ on both $\bamgrpG_2=\Sp_2(q)$ and $\bamgrpG_1=\SL_2(q)$, which is again $\trin$. Clearly these actions correspond to each other via the standard identification maps $\famg_{i,j}$ and $\famg_{j,i}$.

5.~The action of $\widetilde{\Aut}(\FF_{q^2})\le \amgrpC_{i,j}$ on $\amgrpG_i$ and $\amgrpG_j$ via $\ama$ was explained in Remark~\ref{rem:SU2 to SL2 Autos} and Definition~\ref{dfn:check autos of SU4}.
In case $\amgrpG_{i,j}=\SU_4(q)$, $\trin$ is given by conjugation by $M$ composed with the field automorphism 
$\hat{\sigma}$, where $\sigma\colon x\mapsto x^q$ for $x\in \FF_{q^2}$. The same holds for $\bamgrpG_j=\SU_2(q)$ and $\trin$ restricts to $\bamgrpG_i=\SL_2(q^2)$ as transpose-inverse.
In view of Remark~\ref{rem:SU2 to SL2 Autos} we see that via the standard identification map each restricts to 
 transpose inverse on $\amgrpG_i$ and $\amgrpG_j$. 
\epf

\subsection{A standard form for  Curtis-Tits amalgams}\label{subsec:trivial support}
Suppose that $\famG=\{\amgrpG_i,\amgrpG_{i,j},\famg_{i,j}\mid i,j\in I\}$ is a  Curtis-Tits amalgam over $\FF_q$ with $3$-spherical diagram $\liediag$.
Without loss of generality we will assume that all inclusion maps $\famg_{i,j}$ are the standard identification maps of Definition~\ref{dfn:standard CT identification map}.

By Proposition~\ref{prop:existence of weak system of fundamental root groups} it possesses a weak system of fundamental root groups 
\begin{align*}
\Chi=\{\{\amgrpX_i^+,\amgrpX_i^-\}\colon i\in I\},
\end{align*}
 which via the standard embeddings $\famg_{i,j}$ can be identified with those given in Subsection~\ref{subsub:fundamental root groups in CT standard pairs} (note that orienting $\Chi$ may involve changing some signs).
Let $\amA=\{\amgrpA_i,\amgrpA_{i,j},\ama_{i,j}\mid i,j\in I\}$ be the coefficient system associated to $\famG$ and $\Chi$.

We wish to classify all  Curtis-Tits amalgams  $\amG=\{\amgrpG_i,\amgrpG_{i,j},\amg_{i,j}\mid i,j\in I\}$ over $\FF_q$ with the same diagram as $\famG$ with weak system of fundamental root groups $\Chi$ up to isomorphism of  Curtis-Tits amalgams.
By Proposition~\ref{prop:coefficient systems} we may restrict to those amalgams whose connecting maps 
 are of the form $\amg_{i,j}=\famg_{i,j}\after \delta_{i,j}$ for $\delta_{i,j}\in \amgrpA_i$ for all $i\in I$.
 
\bde\label{dfn:trivial support}
The {\dfn trivial support} of $\amG$ (with respect to $\famG$) is the set $\{(i,j)\in I\times I\mid \amg_{i,j}=\famg_{i,j}\}$ (that is, $\delta_{i,j}=\id_{\amgrpG_i}$ in the notation of Proposition~\ref{prop:coefficient systems}). The word ``trivial'' derives from the assumption that the $\famg_{i,j}$'s are the standard identification maps of Definition~\ref{dfn:standard CT identification map}.
\ede
\medskip

Fix some spanning tree $\Sigma\sbe \Gamma$ and suppose that $\edg-\edg\Sigma=\{\{i_s,j_s\}\colon s=1,2,\ldots,r\}$
 so that $H_1(\Gamma,\ZZ)\cong\ZZ^r$. 
 
\bpr\label{prop:trivial support on spanning tree}
There is a  Curtis-Tits amalgam $\amG(\Sigma)$ over $\FF_q$ with the same diagram as $\famG$ and the same $\Chi$, which is isomorphic to $\amG$ and has the following properties:
\begin{enumerate}
\item $\amG$ has trivial support $S=\{(i,j)\in I\times I\mid \{i,j\}\in \edg \Sigma\}\cup \{(i_s,j_s)\colon s=1,2\ldots,r\}$.
\item for each  $s=1,2,\ldots,r$, we have $\amg_{j_s,i_s}=\famg_{j_s,i_s}\after\gamma_{j_s,i_s}$, where $\gamma_{j_s,i_s}\in \amgrpC_{j_s}$.
\end{enumerate}
\epr

\ble\label{lem:stripping tori from connecting maps}
There is a  Curtis-Tits amalgam $\amG^+$ over $\FF_q$ with the same diagram as $\famG$ and the same $\Chi$, which is isomorphic to $\amG$ and has the following properties:
For any $u,v\in I$, if $\amg_{u,v}=\famg_{u,v}\after \gamma_{u,v}\after d_{u,v}$,  for some $\gamma_{u,v}\in \amgrpC_{u}$ and $d_{u,v}\in \amgrpT_u$, then $\amg_{u,v}^+=\famg_{u,v}\after\gamma_{u,v}$.
\ele
\bpf
Note that we have $|I|\ge 2$ and that $\liediag$ is connected.
Fix $u\in I$. Since $\liediag$ is $3$-spherical, there is at most one $w\in I$ such that $(\amgrpG_{u,w},\bamgrpG_u,\bamgrpG_w)$ is a Curtis-Tits standard pair of type $B_2$ or $C_2$.
If there is no such $w$, let $w$ be an arbitrary vertex such that $\{u,v\}\in \edg\liediag$.
We define $\amG^+$ by setting $\amg^+_{u,v}=\famg_{u,v} \after \gamma_{u,v}$  for all $v\ne u$.

Next we define $\phi\colon \amG\to\amG^+$ setting
 $\phi_u=d_{u,w}$ and $\phi_v=\id_{\amgrpG_v}$ for all $v\ne u$.
Now note that setting $\phi_{u,w}=\id_{\amgrpG_{u,v}}$, $\{\phi_{u,w}, \phi_u,\phi_w\}$ is an isomorphism of the subamalgams of $\amG_{\{u,w\}}$ and $\amG^+_{\{u,w\}}$.
As for $\phi_{u,v}$ for $v\ne w$, note that in order for  $\{\phi_{u,v}, \phi_u,\phi_v\}$ to be an isomorphism of the subamalgams of $\amG_{\{u,v\}}$ and $\amG^+_{\{u,v\}}$, we must have 
\begin{align*}
\amg^+_{u,v} \phi_u & = \phi_{u,v}\after\amg_{u,v}\\
\amg^+_{v,u} \phi_v & = \phi_{u,v}\after\amg_{v,u}
\end{align*}
which translates as 
\begin{alignat*}{2}
\famg_{u,v}\after\gamma_{u,v}\after d_{u,w} & = \phi_{u,v}\after\famg_{u,v}\after\gamma_{u,v}\after  d_{u,v}\\
\famg_{v,u} \after\delta_{v,u}& = \phi_{u,v}\after\famg_{v,u}\after\delta_{v,u}
\end{alignat*}
or in other words
\begin{alignat*}{2}
\gamma_{u,v} \after d_{u,w} \after d_{u,v}^{-1}\after \gamma_{u,v}^{-1} & = \ama_{v,u}(\phi_{u,v})\\
\id_{\amgrpG_v} & = \ama_{u,v}(\phi_{u,v}).
\end{alignat*}
Note that $\gamma_{u,v} \after d_{u,w} \after\after d_{u,v}^{-1}\after \gamma_{u,v}^{-1}\in \amgrpT_u\normal \amgrpA_u$. 
Now by Lemma~\ref{lem: coefficient system connecting maps} as 
 $(\amgrpG_{u,v},\bamgrpG_u,\bamgrpG_v)$ is not of type $B_2$ or $C_2$ the map $(\ama_{j,i},\ama_{i,j})\colon \amgrpT_{i,j}\to \amgrpT_i\times\amgrpT_j$ is onto. In particular, the required $\phi_{u,v}\in \amgrpT_{u,v}$ can be found.
This completes the proof.
\epf

\medskip
By Lemma~\ref{lem:stripping tori from connecting maps} in order to prove Proposition~\ref{prop:trivial support on spanning tree} we may now assume that $\amg_{u,v}=\famg_{u,v}\after\gamma_{u,v}$ for some $\gamma_{u,v}\in \amgrpC_u$ for all $u,v\in I$.

Let 
$\amG=\{\amgrpG_i,\amgrpG_j, \amgrpG_{i,j},\amg_{i,j}=\famg_{i,j}\after\gamma_{i,j}\mid i,j\in I\}$ be a Curtis-Tits amalgam over $\FF_q$ with $|I|=2$ and  $\gamma_{i,j} \in \amgrpC_i$ and   $\gamma_{j,i}\in \amgrpC_j$.
We will describe all possible amalgams 
$\amG^+=\{\amgrpG_i,\amgrpG_j, \amgrpG_{i,j},\amg_{i,j}^+=\famg_{i,j}\after\gamma_{i,j}^+\mid i,j\in I\}$ with 
  $\gamma_{i,j}^+ \in \amgrpC_i$ and  $\gamma_{j,i}^+\in \amgrpC_j$, 
isomorphic to $\amG$ via an isomorphism $\phi$ with $\phi_i\in \amgrpC_i$, $\phi_j\in \amgrpC_j$ and 
 $\phi_{i,j}\in \amgrpC_{i,j}$.
 \begin{figure}[H]
\begin{center}
\begin{tikzpicture}
  \matrix[matrix of math nodes,column sep={64pt,between origins},row
%    sep={40pt,between origins},nodes={asymmetrical rectangle}] (s)
    sep={40pt,between origins},nodes={rectangle}] (s)
  {
      &|[name=Gij]|\amgrpG_{i,j} &      \\
 |[name=Gi]|  \amgrpG_i     &  &|[name=Gj]| \amgrpG_j  \\
 |[name=Gpi]|  \amgrpG_i      &  &|[name=Gpj]| \amgrpG_j   \\
     &|[name=Gpij]|\amgrpG_{i,j} &    \\ };
  \draw[ ->]
 % \draw[right hook ->]
   (Gi) edge  node[sloped,above] {$\amg_{i,j}$} (Gij)
                    (Gpj) edge   node[sloped,below] {$\amg^+_{j,i}$}  (Gpij);
%  \draw[left hook ->] 
  \draw[ ->] 
                    (Gj) edge node[sloped,above] {$\amg_{j,i}$} (Gij)
         (Gpi) edge node[sloped,below] {$\amg^+_{i,j}$} (Gpij);
  \draw[->] 
(Gi) edge node[left] {$\phi_i$}  (Gpi)
(Gij) edge [densely dotted] node[left] {$\phi_{i,j}$}  (Gpij)
(Gj) edge node[right] {$\phi_j$}  (Gpj);  
\end{tikzpicture}
\end{center}
\caption{The commuting hexagon of Corollary~\ref{cor:completing the hexagon}.}\label{fig:completing the hexagon}
\end{figure}
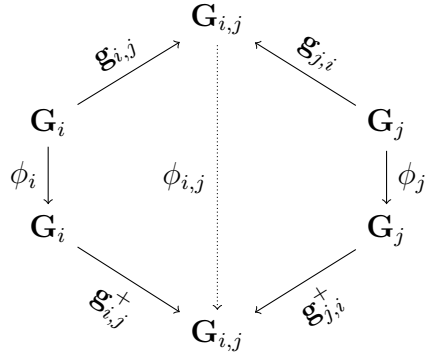
\bco\label{cor:completing the hexagon}

With the notation introduced above, fix the maps $\gamma_{i,j},\gamma^+_{i,j}, \phi_i \in \amgrpC_i$ as well as 
 $\gamma_{j,i}\in \amgrpC_j$. 
Then for any one of $\gamma^+_{j,i}, \phi_j\in \amgrpC_j$, there exists a choice $\gamma\in \amgrpC_i$ for the remaining map in $\amgrpC_j$  so that there exists $\phi_{i,j}$ making the diagram in Figure~\ref{fig:completing the hexagon} commute.
Moreover, if $\liediag_{i,j}$ is one of $A_2$, $B_2$, $C_2$, $\twA_3$, then $\gamma$ is unique, whereas if $\liediag_{i,j}=\twD_3$, then there are exactly two choices for $\gamma$. 
\eco
\bpf
The first claim follows immediately from the fact that the restriction maps $\ama_{j,i}\colon \amgrpC_{i,j}\to \amgrpC_i$ and $\ama_{i,j}\colon \amgrpC_{i,j}\to \amgrpC_j$ in part 4.~and 5.~of Lemma~\ref{lem: coefficient system connecting maps}  are both surjective. The second claim follows from the fact that $\ama_{j,i}\colon \amgrpC_{i,j}\to\amgrpC_i$ is injective except if $\Gamma_{i,j}=\twD_3$ in which case it has a kernel of order $2$.
\epf

\medskip

\bpf (of Proposition~\ref{prop:trivial support on spanning tree})
By Lemma~\ref{lem:stripping tori from connecting maps} we may assume that $\amg_{i,j}=\famg_{i,j}\after\gamma_{i,j}$ for some $\gamma_{i,j}\in \amgrpC_i$ for all $i,j\in I$.

For any (possibly empty) subset $T\sbe\vrtc$ let $S(T)$ be the set of pairs $(i,j)\in S$ such that $i\in T$.
Clearly the trivial support of $\amG$ contains $S(\emptyset)$.

We now show that if $T$ is the vertex set of a (possibly empty) proper subtree of $\Sigma$, and $u$ is a vertex such that $T\cup \{u\}$ is also the vertex set of a subtree of $\Sigma$, then for any  Curtis-Tits amalgam $\amG$ whose trivial support contains $S(T)$, there is a  Curtis-Tits amalgam $\amG^+$ isomorphic to $\amG$, whose trivial support contains $S(T\cup\{u\})$.

Once this is proved, Claim 1.~follows since we can start with $T=\emptyset$ and end with a  Curtis-Tits amalgam, still isomorphic to $\amG$, whose trivial support contains $S$.

Now let $T$ and $u$ be as above. 
We first deal with the case where $T\ne\emptyset$. 
Let $t$ be the unique neighbor of $u$ in the subtree of $\Sigma$ with vertex set $T\cup\{u\}$.
We shall define an amalgam $\amG^+=\{\amgrpG_i,\amgrpG_{i,j},\amg_{i,j}^+=\famg_{i,j}\after\gamma^+_{i,j}\mid i,j\in I\}$ and an isomorphism
 $\phi\colon \amG\to \amG^+$, where $\gamma^+_{i,j},\phi_i\in \amgrpC_i$ and $\phi_{\{i,j\}}\in \amgrpC_{i,j}$ for all $i,j\in I$. 
First note that it suffices to define $\amg^+_{i,j}$, $\phi_i$ and $\phi_{\{i,j\}}$ for $\{i,j\}\in \edg$: given this data, by the $A_1\times A_1$ case in Lemma~\ref{lem:coefficient system groups} and Lemma~\ref{cor:completing the hexagon}, for any non-edge $\{k,l\}$ there is a unique $\phi_{\{k,l\}}\in \amgrpC_{k,l}$ such that 
 $(\phi_{k,l},\phi_k,\phi_l)$ is an isomorphism between  $\amG_{\{k,l\}}$ and $\amG^+_{\{k,l\}}$.

Before defining inclusion maps on edges, note that since $\liediag$ is $3$-spherical, no two neighbors of $u$ in $\Gamma$ are connected by an edge.
Therefore we can unambiguously set 
\begin{align*}
\amg^+_{i,j} & = \amg_{i,j}  \mbox{ for }u\not\in\{i,j\}\in \edg\liediag. 
\end{align*} 
Note that both maps $\amg^+_{t,u}$ and $\amg^+_{u,t}$ are forced upon us, but at  this point for any other neighbor $v$ of $u$, only one of $\amg^+_{u,v}$ and $\amg^+_{v,u}$ is forced upon us. 
We set 
\begin{align*}
\amg^+_{t,u}&=\amg_{t,u}, \mbox{ and }\\
\amg^+_{v,u}&=\amg_{v,u} \mbox{ for }v\in I \mbox{ with }(u,v)\not\in S \mbox{ and }(v,u)\in S.
\end{align*}
To extend the trivial support as required, we set 
\begin{align*}
\amg^+_{u,v}&=\famg_{u,v} \mbox{ for }v\in I \mbox{ with }(u,v)\in S.
\end{align*}
We can already specify part of $\phi$: Set
\begin{align*}
\phi_i& =\id_{\amgrpG_i} \mbox{ for }i\in I-\{u\},\\
\phi_{\{i,j\}}&=\id_{\amgrpG_{i,j}} \mbox{ for } u\not\in \{i,j\}\in \edg\liediag.
\end{align*}
Thus, what is left to specify is the following: $\phi_u$ and $\phi_{\{u,t\}}$ and, for all neighbors $v\ne t$ of $u$ we must specify $\phi_{\{u,v\}}$ as well as  
\begin{align*}
\amg^+_{u,v} & \mbox{ if  } (u,v)\not\in S,\\
\amg^+_{v,u} & \mbox{ if }  (u,v)\in S.
\end{align*}  

Figures~\ref{fig:vu not in S}~and~\ref{fig:uv not in S} describe the amalgam $\amG$ (top half) and $\amG^+$ (bottom half) at the vertex $u$, where $t\in \vrtc \Sigma$, $v\in \vrtc \liediag$, and $\{u,t\}, \{u,v\}\in \edg\liediag$.
Inclusion maps from $\amG^+$ forced upon us are indicated in bold, the dotted arrows are those we must define so as to make the diagram commute.

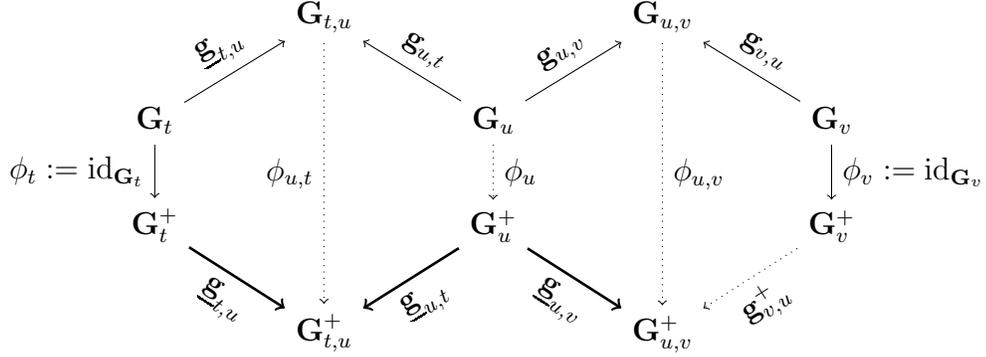
\begin{figure}[H]
\begin{center}
\begin{tikzpicture}
  \matrix[matrix of math nodes,column sep={64pt,between origins},row
%    sep={40pt,between origins},nodes={asymmetrical rectangle}] (s)
    sep={40pt,between origins},nodes={rectangle}] (s)
  {
      &|[name=Gtu]|\amgrpG_{t,u} & &|[name=Guv]| \amgrpG_{u,v} &     \\
 |[name=Gt]|  \amgrpG_t      &  &|[name=Gu]| \amgrpG_u& & |[name=Gv]| \amgrpG_v    \\
 |[name=Gpt]|  \amgrpG^+_t      &  &|[name=Gpu]| \amgrpG^+_u&  & |[name=Gpv]| \amgrpG^+_v    \\
     &|[name=Gptu]|\amgrpG^+_{t,u} & &|[name=Gpuv]| \amgrpG^+_{u,v} &     \\
  };
%  \draw[right hook ->]
  \draw[ ->]
   (Gt) edge  node[sloped,above] {$\famg_{t,u}$} (Gtu)
              (Gu) edge  node[sloped,above] {$\amg_{u,v}$} (Guv)
                    (Gpu) edge [line width = 1pt]  node[sloped,below] {$\famg_{u,t}$}  (Gptu)
   (Gpv) edge [dotted] node[sloped,below] {$\amg^+_{v,u}$} (Gpuv)
 ;
 
 %\draw[left hook ->] 
 \draw[->] 
                    (Gu) edge   node[sloped,above] {$\amg_{u,t}$} (Gtu)
     (Gv) edge  node[sloped,above] {$\amg_{v,u}$} (Guv)
         (Gpu) edge  [line width = 1pt]  node[sloped,below] {$\famg_{u,v}$}  (Gpuv)
         (Gpt) edge  [line width = 1pt]  node[sloped,below] {$\famg_{t,u}$} (Gptu)
    ;
  
 \draw[->] 
(Gt) edge  node[left] {$\phi_t:=\id_{\amgrpG_t}$}  (Gpt)
(Gtu) edge [ dotted]  node[left] {$\phi_{u,t}$}  (Gptu)
(Gu) edge [ dotted]  node[right] {$\phi_u$}  (Gpu)
(Guv) edge [ dotted]  node[right] {$\phi_{u,v}$}  (Gpuv)
(Gv) edge  node[right] {$\phi_v:=\id_{\amgrpG_v}$}  (Gpv)
;  
\end{tikzpicture}
\caption{The case $(u,v)\in S$ and $(v,u)\not\in S$.}\label{fig:vu not in S}
\end{center}
\end{figure}
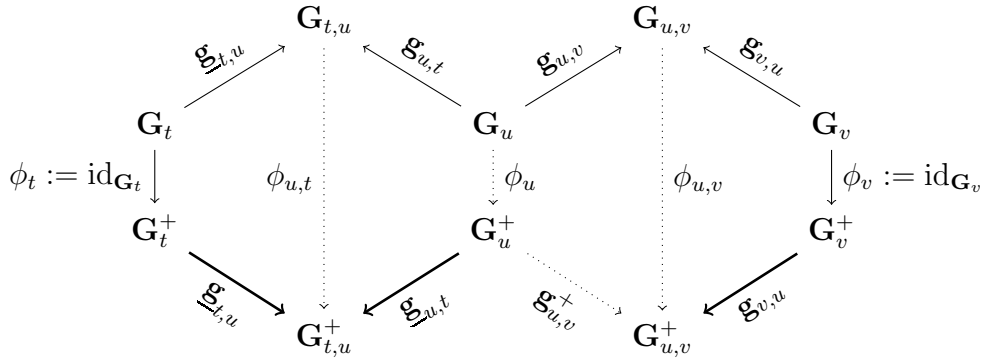
\begin{figure}[H]
\begin{center}
\begin{tikzpicture}
  \matrix[matrix of math nodes,column sep={64pt,between origins},row
%    sep={40pt,between origins},nodes={asymmetrical rectangle}] (s)
    sep={40pt,between origins},nodes={rectangle}] (s)
  {
      &|[name=Gtu]|\amgrpG_{t,u} & &|[name=Guv]| \amgrpG_{u,v} &     \\
 |[name=Gt]|  \amgrpG_t      &  &|[name=Gu]| \amgrpG_u& & |[name=Gv]| \amgrpG_v    \\
 |[name=Gpt]|  \amgrpG^+_t      &  &|[name=Gpu]| \amgrpG^+_u&  & |[name=Gpv]| \amgrpG^+_v    \\
     &|[name=Gptu]|\amgrpG^+_{t,u} & &|[name=Gpuv]| \amgrpG^+_{u,v} &     \\
  };
%  \draw[right hook ->]
  \draw[ ->]
   (Gt) edge  node[sloped,above] {$\famg_{t,u}$} (Gtu)
              (Gu) edge  node[sloped,above] {$\amg_{u,v}$} (Guv)
                    (Gpu) edge [line width = 1pt]  node[sloped,below] {$\famg_{u,t}$}  (Gptu)
   (Gpv) edge [line width = 1pt] node[sloped,below] {$\amg_{v,u}$} (Gpuv)
 ;
 
% \draw[left hook ->] 
 \draw[ ->] 
                    (Gu) edge   node[sloped,above] {$\amg_{u,t}$} (Gtu)
     (Gv) edge  node[sloped,above] {$\amg_{v,u}$} (Guv)
         (Gpu) edge [dotted] node[sloped,below] {$\amg^+_{u,v}$}  (Gpuv)
         (Gpt) edge  [line width = 1pt]  node[sloped,below] {$\famg_{t,u}$} (Gptu)
    ;
  
 \draw[->] 
(Gt) edge  node[left] {$\phi_t:=\id_{\amgrpG_t}$}  (Gpt)
(Gtu) edge [ dotted]  node[left] {$\phi_{u,t}$}  (Gptu)
(Gu) edge [ dotted]  node[right] {$\phi_u$}  (Gpu)
(Guv) edge [ dotted]  node[right] {$\phi_{u,v}$}  (Gpuv)
(Gv) edge  node[right] {$\phi_v:=\id_{\amgrpG_v}$}  (Gpv)
;  
\end{tikzpicture}
\caption{The case $(u,v)\not\in S$ and $(v,u)\in S$.}\label{fig:uv not in S}
\end{center}
\end{figure}
In these figures all non-dotted maps are of the form $\famg_{i,j}\after\gamma_{i,j}$ for some $\gamma_{i,j}\in \amgrpC_i$ hence we can find the desired maps using Corollary~\ref{cor:completing the hexagon}.

In case $T=\emptyset$, the situation is as described in Figures~\ref{fig:vu not in S}~and~\ref{fig:uv not in S} after removing the $\{u,t\}$-hexagon and any conditions it may impose on $\phi_u$, and letting $v$ run over all neighbors of $u$. That is, we must now define $\phi_u$, and for any neighbor $v$ of $u$, we must find $\phi_{u,v}$ as well as 
\begin{align*}
\amg^+_{u,v} & \mbox{ if  } (u,v)\not\in S,\\
\amg^+_{v,u} & \mbox{ if }  (u,v)\in S.
\end{align*}
To do so we let $\phi_u=\id_{\amgrpG_u}\in \amgrpC_u$. 
Finally, for each neighbor $v$ of $u$ we simply let $\amg_{u,v}^+=\amg_{u,v}$ (so that $\phi_{u,v}=\id_{\amgrpG_{i,j}}\in \amgrpC_{i,j}$) if $(u,v)\not\in S$, and we obtain $\amg^+_{v,u}$ and $\phi_{u,v}\in \amgrpC_{u,v}$ using Corollary~\ref{cor:completing the hexagon} if $(u,v)\in S$.
\epf

\subsection{Classification of  Curtis-Tits amalgams with $3$-spherical diagram}\label{subsec:classification of 3-spherical CT amalgams}

In the case where $\amG$ is a Curtis-Tits amalgam over $\FF_q$ whose diagram is a $3$-spherical tree,  Proposition~\ref{prop:trivial support on spanning tree} says that $\amG\cong\famG$. 

\bth\label{thm:CT 3 spherical tree}
Suppose that $\amG$ is a  Curtis-Tits amalgam with a diagram that is a $3$-spherical tree.
Then, $\amG$ is unique up to isomorphism. 
In particular any  Curtis-Tits amalgam with spherical diagram is unique.
\eth

\ble\label{lem:spanning tree with minimal complement}
Given a Curtis-Tits amalgam over $\FF_q$ with connected $3$-spherical diagram $\liediag$ there is a spanning tree $\Sigma$ such that the set of edges 
 in $\edg\liediag-\edg\Sigma=\{\{i_s,j_s\}\colon s=1,2,\ldots,r\}$ has the property that
 \begin{enumerate}
 \item\label{cond:A2} $(\amgrpG_{\{i_s,j_s\}},\amg_{i_s,j_s}(\amgrpG_{i_s}),\amg_{i_s,j_s}(\amgrpG_{j_s}))$ has type $A_2(q^{e_s})$, where $e_s$ is some power of $2$.
  \item\label{cond:minimal e} There is a loop $\Lambda_s$ containing $\{i_s,j_s\}$ such that 
  any vertex group of $\Lambda_s$ is isomorphic to $\SL_2(q^{e_s 2^l})$ for some $l\ge 0$.  
 \end{enumerate}
\ele
\bpf
Induction on the rank $r$ of  $H^1(\liediag,\ZZ)$.
If $r=0$, then there is no loop at all and we are done.

Consider the collection of all edges $\{i,j\}$ of $\Gamma$ such that $\liediag_{\{i,j\}}$ has type $A_2$ and 
 $H^1(\liediag-\{i,j\},\ZZ)$ has rank $r-1$, and choose one such that $\amgrpG_i\cong\SL_2(q^{e_1})$ where $e_1$ is minimal among all these edges.  
Next replace $\liediag$ by $\liediag-\{i,j\}$ and use induction.
Suppose $\{i_s,j_s\}\mid s=1,2,\ldots,r\}$ is the resulting selection of edges so that $\Sigma=\liediag-\{i_s,j_s\}\mid s=1,2,\ldots,r\}$ is a spanning tree and condition~\eqref{cond:A2} is satisfied.
Note that by choice of these edges, also condition~\eqref{cond:minimal e} is satisfied by at least one of the loops of $\liediag-\{\{i_t,j_t\}\colon t=1,2,\ldots,s-1\}$ that contains $\{i_s,j_s\}$.
Note that this uses the fact that by $3$-sphericity every vertex belongs to at least one subdiagram of type $A_2$.
\epf

\bde\label{dfn:CT the map kappa}
Fix a connected $3$-spherical diagram $\liediag$ and a prime power $q$. Let  $\Sigma$  be a spanning tree and let the set of edges $\edg\liediag-\edg\Sigma=\{\{i_s,j_s\}\colon s=1,2,\ldots,r\}$ together with the integers $\{e_s\colon s=1,2,\ldots,r\}$ satisfy the conclusions of Lemma~\ref{lem:spanning tree with minimal complement}.
Let $\sCT(\liediag,q)$ be the collection of isomorphism classes of  Curtis-Tits amalgams of type $\liediag(q)$ and let $\famG=\{\amgrpG_i,\amgrpG_{i,j},\famg_{i,j}\mid i,j\in I\}$ be the standard Curtis-Tits amalgam over $\FF_q$ with diagram $\liediag$ as in Subsection~\ref{subsec:trivial support}.

Consider the following map:
\begin{align*}
\kappa\colon\prod_{s=1}^r \Aut(\FF_{q^{e_s}})\times \langle\trin\rangle\to \sCT(\liediag).
\end{align*}
where $\kappa((\alpha_s)_{s=1}^r)$ is the isomorphism class of the amalgam $\amG^+=\amG((\alpha_s)_{s=1}^r)$ given by setting $\amg^+_{j_s,i_s}=\famg_{j_s,i_s}\after \alpha_s$ for all $s=1,2,\ldots,r$.
\ede

We now have 
\bco\label{cor:CT kappa is onto}
The map $\kappa$ is onto.
\eco
\bpf
Note that, for each $s=1,2,\ldots,r$, the Curtis-Tits standard pair $(\amgrpG_{\{i_s,j_s\}}$, $\amg_{i_s,j_s}(\amgrpG_{i_s})$, $\amg_{i_s,j_s}(\amgrpG_{j_s}))$ has type $A_2(q^{e_s})$ and so $\amgrpC_{j_s}=\Aut(\FF_{q^{e_s}})$.  Thus the claim is an immediate consequence of Proposition~\ref{prop:trivial support on spanning tree}.
\epf

\medskip
We note that if we select $\Sigma$ differently, the map $\kappa$ will still be onto. However, the ``minimal'' choice made in Lemma~\ref{lem:spanning tree with minimal complement} ensures that $\kappa$ is injective as well, as we will see.

\ble\label{lem:classification of  CT amalgams with loop diagram}
Suppose $\liediag(q)$ is a 
$3$-spherical diagram $\liediag$ that is a simple loop. 
Then, $\kappa$ is injective.
\ele
\bpf
Suppose there is an isomorphism $\kappa(\alpha)=\amG\stackrel{\phi}{\longrightarrow} \amG^+=\kappa(\beta)$, for some  $\alpha,\beta\in \Aut(\FF_q)\times\langle \trin\rangle$.
 Write $I=\{0,1,\ldots,n-1\}$ so that $\{i,i+1\}\in \edg\liediag$ for all $i\in I$ (subscripts modulo $n$).
 Without loss of generality assume that $(i_1,j_1)=(1,0)$ so that by Proposition~\ref{prop:trivial support on spanning tree} we  may assume that $\amg_{i,j}=\famg_{i,j}=\amg^+_{i,j}$ for all $(i,j)\ne (1,0)$.
This means that  $\ama\colon \amgrpC_{i,i+1}\to \amgrpC_i\times\amgrpC_{i+1}$ sends $\phi_{i,i+1}$ to $(\phi_i,\phi_{i+1})$ for any edge $\{i,i+1\}\ne \{0,1\}$. 
Now note that by minimality of $q$, $\amgrpC_i$  (and $\amgrpC_{i,i+1}$) has a quotient $\bar{\amgrpC}_i$   (and $\bar{\amgrpC}_{i,i+1}$) isomorphic to $\langle\Aut(\FF_q)\rangle\times\langle\trin\rangle$ for every $i\in I$, by considering the action of $\amgrpC_i$ on the subgroup of $\amgrpG_i$ isomorphic to $\SL_2(q)$. 
By Part 4 and 5 of Lemma~\ref{lem: coefficient system connecting maps} the maps  $\ama_{i+1,i}^{-1}$ and $\ama_{i,i+1}$ induce isomorphisms $\bar{\amgrpC}_i\to \bar{\amgrpC}_{i,i+1}$ and $\bar{\amgrpC}_{i,i+1}\to\bar{\amgrpC}_{i+1}$, which compose to an isomorphism
\begin{align*}
\phi_i&\mapsto \famg_{i+1,i}^{-1}\after \famg_{i,i+1}\after \phi_i \after \famg_{i,i+1}^{-1}\after \famg_{i+1,i},
\end{align*}
sending the image of $\trin$ and $\alpha$ in $\bar{\amgrpC}_i$ to the image of $\trin$ (and $\alpha$ respectively) in $\bar{\amgrpC}_{i+1}$, where $\alpha\colon x\mapsto x^p$ for $x$ in the appropriate extension of $\FF_q$ defining $\amgrpG_{i,i+1}$. 
Concatenating these maps along the path $\{1,2,\ldots,n-1,0\}$ and considering the edge $\{0,1\}$ we see that the images of $\beta^{-1}\phi_1\alpha$ and $\phi_1$ in $\bar{\amgrpC}_1=\amgrpC_1$ coincide.
 Since $\amgrpC_1$ is abelian this means that $\beta=\alpha$.
\epf

\bth\label{thm:CT classification of 3-spherical amalgams}
Let  $\liediag$ be a connected $3$-spherical diagram with spanning tree $\Sigma$ and set of edges
 $\edg\liediag-\edg\Sigma=\{\{i_s,j_s\}\colon s=1,2,\ldots,r\}$ together with the integers $\{e_s\colon s=1,2,\ldots,r\}$ satisfying the conclusions of Lemma~\ref{lem:spanning tree with minimal complement}.
Then $\kappa$ is a bijection between the elements of $\prod_{s=1}^r \Aut(\FF_{q^{e_s}})\times \langle\trin\rangle$ and the type preserving isomorphism classes of  Curtis-Tits amalgams with diagram $\liediag$ over $\FF_q$.
\eth
\bpf
Again, it suffices to show that $\kappa$ is injective. This in turn follows from Lemma~\ref{lem:classification of  CT amalgams with loop diagram}, for if two amalgams are isomorphic (via a type preserving isomorphism), then the amalgams induced on subgraphs of $\liediag$ must be isomorphic and Lemma~\ref{lem:classification of  CT amalgams with loop diagram} shows that $\kappa$ is injective on the subamalgams supported by the loops $\Lambda_s$ ($s=1,2,\ldots,r$).
\epf
\section{Classification of  Phan amalgams}\label{sec:classification of Ph amalgams}
\subsection{Introduction}
The classification problem is formulated as follows: Determine, up to isomorphism of amalgams, all  Phan amalgams $\amG$ with given diagram $\liediag$ possessing a non-trivial (universal) completion.

\subsection{Classification of Phan amalgams with $3$-spherical diagram}
\subsubsection{Tori in Phan standard pairs}\label{subsubsec:Tori Phan standard pairs}
Let $\amG=\{\amgrpG_{i,j},\amgrpG_i,\amg_{i,j}\mid i,j\in I\}$ be a Phan amalgam over $\FF_q$ with $3$-spherical diagram $\liediag=(I,E)$.
This means that the subdiagram of $\liediag$ induced on any set of three vertices is spherical.
This is equivalent to $\liediag$ not containing triangles of any kind and such that no vertex is on more than one $C_2$-edge.

\bde\label{dfn:tori in phan amalgam}
For any $i,j\in I$ with $\{i,j\}\in \edg\liediag$, let 
\begin{align*}
\amgrpD_i^j=N_{\amgrpG_{i,j}}(\amg_{j,i}(\amgrpG_j))\cap \amg_{i,j}(\amgrpG_i)
\end{align*}
\ede
\ble\label{lem:tori in phan standard pairs}
Suppose that $(\amgrpG,\amgrpG_1,\amgrpG_2)$ is a Phan standard pair of type 
$\liediag(q)$ as in Subsection~\ref{subsec:standard P pairs}. 
\begin{enumerate}
\item If $\liediag(q)=A_2(q)$, then $\langle \amgrpD_1^2, \amgrpD_2^1\rangle$ is the standard torus stabilizing the orthonormal basis $\{e_1,e_2,e_3\}$. Here $\amgrpD_1^2$ (resp. $\amgrpD_2^1$) is the stabilizer in this torus of $e_1$ (resp. $e_3$).
\item If $\liediag(q)=C_2(q)$, then $\langle \amgrpD_1^2,\amgrpD_2^1\rangle$ is the standard torus stabilizing the 
 basis $\{e_1,e_2,e_3=f_1,e_4=f_2\}$ which is hyperbolic for the symplectic form of $\Sp_4(q^2)$ and orthonormal for the unitary form of $\SU_4(q)$; Here $\amgrpD_2^1$ (resp. $\amgrpD_1^2$) is the stabilizer of $\langle e_2\rangle$ and $\langle f_2\rangle$ (resp.~the pointwise stabilizer of both $\langle e_1,f_1\rangle$ and $\langle e_2,f_2\rangle$).
Thus, 
\begin{align*}
\amgrpD_2^1&=\langle\diag(1,a,1,a^\sigma)\colon a\in \FF_{q^2}\mbox{ with }aa^\sigma=a^{q+1}=1\rangle,\\
\amgrpD_1^2&=\langle\diag(a,a^\sigma,a^\sigma,a )\colon a\in \FF_{q^2}\mbox{ with }aa^\sigma=a^{q+1}=1\rangle.
\end{align*}
\item In either case, for $\{i,j\}=\{1,2\}$,  $\amgrpD_i^j=C_{\amgrpG_{i,j}}(\amgrpD_j^i)\cap\amgrpG_i$ and 
 $\amgrpD_j^i$ is the unique torus of $\amgrpG_j$ normalized by $\amgrpD_i^j$. 
\end{enumerate}
\ele
\bpf
Parts 1.~and~2.~as well as the first claim of Part 3.~are straightforward matrix calculations.
As for the last claim note that in both cases, $\amgrpD_i^j$ acts diagonally on $\amgrpG_j$ viewed as $\SU_2(q)$ in it natural representation $V$ via the standard identification map; in fact  (in the case $C_2$, $\amgrpD_2^1$ acts even innerly on $\amgrpG_1$).
If $\amgrpD_i^j$ normalizes a torus $\amgrpD'$ in $\amgrpG_j$ then it will have to stabilize its eigenspaces.
Since $q+1\ge 3$, the eigenspaces of $\amgrpD_i^j$ in its action on $V$ have dimension $1$, so $\amgrpD_i^j$ and $\amgrpD'$ must share these eigenspaces. This means that $\amgrpD'=\amgrpD_j^i$.
\epf

\subsubsection{Property (D) for Phan amalgams}\label{subsubsec:Phan property D}
We state Property (D) for $3$-spherical Phan amalgams, extending the definition from~\cite{BloHof2014b} which was given for Curtis-Tits amalgams under the assumption that $q\ge 4$.
\bde\label{dfn:property D}\nom{{(\rm property (D))}}{}
We say that $\amG$ {\dfn has property (D)} if there is a {\dfn system of tori} $\cD=\{\amgrpD_i\colon i\in I\}$ such that for all edges $\{i,j\}\in \edg\liediag$ we have  $\amg_{i,j}(\amgrpD_i)=\amgrpD_i^j$.
\ede
\ble\label{lem:non-collapsing Phan has property D}
Suppose that $\amG$ has a completion $(\compG,\compg)$ so that $\compg_i$ is non-trivial for all $i\in I$.
Then, for any $i,j,k\in I$ such that $\{i,j\},\{j,k\}\in \edg\liediag$, there is a torus $\amgrpD_j\le \amgrpG_j$ such that 
 $\amg_{j,i}(\amgrpD_j)=\amgrpD_j^i$ and $\amg_{j,k}(\amgrpD_j)=\amgrpD_j^k$.
In particular, $\amG$ has property (D).
\ele
\bpf
First note that in case $q=2$, the conclusion of the lemma is trivially true as, for all $i\in I$, $\amgrpG_i\cong S_3$ has a unique Phan torus.

We now consider the general case.
For $\liediag(q)=A_3(q)$, this was proved by Bennett and Shpectorov in~\cite{BeSh2004} (see also~\cite{BloHof2014b}).
For completeness we recall the argument, which applies in this more general case as well.
We shall prove that 
 \begin{align*}
 \compg(\amgrpD_j^i)&=\compg(\amgrpD_j^k)
 \end{align*}
and then let  $\amgrpD_j\le \amgrpG_j$ be such that $\compg(\amgrpD_j)= \compg(\amgrpD_j^i)=\compg(\amgrpD_j^k)$.
Note that since $\compg_j$ is non-trivial, it now follows that $\amg_{j,i}(\amgrpD_j)=\amgrpD_j^i$ and $\amg_{j,k}(\amgrpD_j)=\amgrpD_j^k$.

Recall that for any subgroup $\amgrpH$ of a group in $\amG$ we'll write $H=\compg(\amgrpH)$.
We show that $D_j^i$ is normalized by $D_k^j$ and use Lemma~\ref{lem:tori in phan standard pairs} to conclude that 
 $D_j^i=D_j^k$.
To that end we let $h\in D_k^j$  and prove that $hD_j^ih^{-1}=D_j^i$. To achieve this we show that $hD_j^ih^{-1}$ is normalized by $D_i^j$ and again use Lemma~\ref{lem:tori in phan standard pairs}.
So now let $g\in D_i^j$ and 
note that since $\liediag$ is $3$-spherical, $\{i,k\}\not\in \edg\liediag$ so that $g$ and $h$ commute.
In addition note that  by Lemma~\ref{lem:tori in phan standard pairs},  $gD_j^ig^{-1}=D_j^i$.
Therefore we have 
\begin{align*}
ghD_j^i h^{-1}g^{-1}= hg D_j^i g^{-1}h^{-1}= h D_j^i h^{-1},
\end{align*}
as required.
\epf

\subsubsection{The coefficient system of a  Phan amalgam}\label{subsec:Phan coefficient system}

\bde\label{dfn:type of P amalgam}
We now fix a standard Phan amalgam  $\famG=\{\amgrpG_i,\amgrpG_{i,j},\famg_{i,j}\mid i,j\in I\}$ over $\FF_q$ with diagram $\liediag(q)$, where for every $i,j\in I$, $\famg_{i,j}$ is the standard identification map of Definition~\ref{dfn:standard Phan identification map}.
Then,  $\famG$ has property (D) with system of tori $\cD=\{\amgrpD_i\colon i\in I\}$ as in Lemma~\ref{lem:tori in phan standard pairs}. 

If $\amG$ is any other  non-collapsing Phan amalgam over $\FF_q$ with diagram $\liediag$, then since all tori of $\amgrpG_i$ are conjugate under $\Aut(\amgrpG_i)$, by adjusting the inclusion maps $\amg_{i,j}$ we can replace $\amG$ by an isomorphic amalgam whose system of tori is exactly $\cD$.
\ede

From now on we assume that $\famG$, $\cD=\{\amgrpD_i\colon i\in I\}$ and $\amG$ are as in Definition~\ref{dfn:type of P amalgam}

\bde\label{dfn:phan coefficient system}
Suppose that $\amG=\{\amgrpG_i,\amgrpG_{i,j},\amg_{i,j}\mid i,j\in I\}$ is a Phan amalgam with connected $3$-spherical diagram $\liediag$ having property (D).
Let $\cD=\{\amgrpD_i\colon i\in I\}$  be the associated system of  tori.
The {\dfn coefficient system associated to $\amG$} is the collection 
 $\amA=\{\amgrpA_i,\amgrpA_{i,j},\ama_{i,j}\mid i,j\in I\}$ where, for any $i,j\in I$ we set 
 \begin{align*}
 \amgrpA_i&=N_{\Aut(\amgrpG_i)}(\amgrpD_i), \\
 \amgrpA_{i,j}&=N_{\Aut(\amgrpG_{i,j})}(\amg_{i,j}(\amgrpG_i))\cap N_{\Aut(\amgrpG_{i,j})}(\amg_{j,i}(\amgrpG_j)), \\
 \ama_{i,j}&\colon \amgrpA_{i,j}\to\amgrpA_j \mbox{ is given by restriction: } \varphi\mapsto \amg_{j,i}^{-1}\after \rho_{i,j}(\varphi)\after\amg_{j,i}.
  \end{align*}
where $\rho_{i,j}(\varphi)$ is the restriction of $\varphi$ to $\bamgrpG_j\le \amgrpG_{i,j}$.
\ede
From now on we let $\amA$ be the coefficient system associated to $\famG$ with respect to the system of tori $\cD$.
The fact that the $\ama_{i,j}$ are well-defined follows from the following simple observation.
\ble\label{lem:N G1 G2=N D1 D2} 
For any $i,j\in I$ with $\{i,j\}\in \edg\liediag$, we have 
\begin{align*}
\amgrpA_{i,j}\le N_{\Aut(\amgrpG_{i,j})}(\amg_{i,j}(\amgrpD_i))\cap N_{\Aut(\amgrpG_{i,j})}(\amg_{j,i}(\amgrpD_j)).
\end{align*}
\ele
\bpf
The inclusion $\le$ is immediate from the definitions.
\epf

\medskip

The significance for the classification of  Phan  amalgams with the same system of tori is as follows:
\bpr\label{prop:P coefficient systems}
Suppose that $\amG$ and $\amG^+$ are  Phan amalgams of type $\famG$ with the same system of tori $\cD=\{\amgrpD_i\colon i\in I\}$.
\begin{enumerate}
\item For all $i,j\in I$, we have $\amg_{i,j}=\famg_{i,j}\after\delta_{i,j}$ and $\amg_{i,j}=\famg_{i,j}\after\delta_{i,j}^+$ for some $\delta_{i,j}, \delta_{i,j}^+\in \amgrpA_i$,
\item For any isomorphism $\phi\colon \amG\to\amG^+$ and $i,j\in I$, we have $\phi_i\in \amgrpA_i$, $\phi_{\{i,j\}}\in \amgrpA_{i,j}$, and $\ama_{i,j}(\phi_{\{i,j\}})=\delta_{i,j}^+\after\phi_i\after\delta_{i,j}^{-1}$.
\end{enumerate}
\epr
\bpf
Part 1.~follows since, for any $i,j\in I$ we have $\amg_{i,j}^{-1}\after\famg_{i,j}\in \Aut(\amgrpG_i)$ and  
\begin{align*}
\amg_{i,j}(\amgrpD_i)=\famg_{i,j}(\amgrpD_i).
\end{align*}
Part 2.~follows from Lemma~\ref{lem:N G1 G2=N D1 D2} since, for any $i,j\in I$,  
\begin{align*}
(\amgrpG_{i,j},\amg_{i,j}(\amgrpG_i),\amg_{j,i}(\amgrpG_j))
=(\amgrpG_{i,j},\famg_{i,j}(\amgrpG_i),\famg_{j,i}(\amgrpG_j)) =(\amgrpG_{i,j},\amg^+_{i,j}(\amgrpG_i),\amg^+_{j,i}(\amgrpG_j)).
\end{align*}
\epf

We now determine the groups appearing in a coefficient system by looking at standard pairs.
\ble\label{lem:structure of Phan A_i groups}
Fix $i\in I$ and let $q$ be such that $\amgrpG_i\cong\SU_2(q)$.
Then, 
\begin{align*}
\amgrpA_i&=\amgrpT_i\rtimes \amgrpC_i,
\end{align*}
where $\amgrpT_i$ is the subgroup of diagonal automorphisms in $\PGU_2(q)$ and $\amgrpC_i=\Aut(\FF_{q^2})$.
%\langle \trin, \Aut(\FF_q)\rangle$.
\ele
\bpf
This follows from the fact that via the standard embedding map $\famg_{i,j}$ the groups $\amgrpD_i$  of the system of tori are the subgroups of standard diagonal matrices in $\SU_2(q)$.

To see this note that $\amgrpG_i\cong\SU_2(q)$ and that $\Aut(\amgrpG_i)\cong\PGU_2(q)\rtimes\Aut(\FF_{q^2})$.
Also, $\amgrpD_i=\langle d\rangle$ for some $d=\diag(\zeta,\zeta^q)$ and $\zeta$ a primitive $q+1$-th root of $1$ in $\FF_{q^2}$.
A quick calculation now shows that $\tau$ and  $\sigma$ are the same in their action, which is inner and one verifies that $N_{\GU_2(q)}(\amgrpD_i)=\langle \trin,\diag(a,b)\colon a,b\in \FF_q^2\rangle$.
\epf

\ble\label{lem:structure of Phan A_ij groups}
Let $\amA$ be the coefficient system associated to the standard Phan amalgam $\famG$ of type $\liediag(q)$ and the system of tori $\cD$.

If $\Gamma=A_1\times A_1$, we have $\amgrpG_{i,j}=\amgrpG_i\times\amgrpG_j$, $\famg_{i,j}$ and $\famg_{j,i}$ are identity maps, and  
\begin{align}
\amgrpA_{i,j}&=\amgrpA_i\times\amgrpA_j\cong \amgrpT_{i,j}\rtimes\amgrpC_{i,j}.\label{eqn:Phan N Xi Xj A1timesA1}
\end{align}
where $\amgrpT_{i,j}=\amgrpT_i\times\amgrpT_j$ and $\amgrpC_{i,j}=\amgrpC_i\times\amgrpC_j$.
Otherwise, 
\begin{align*}
\amgrpA_{i,j}&=\amgrpT_{i,j}\rtimes\amgrpC_{i,j},
\end{align*}
where 
 $\amgrpC_{i,j}=\Aut(\FF_{q^2})$ and  $\amgrpT_{i,j}$ denotes the image of the standard torus $\GD$ in $\Aut(\amgrpG_{i,j})$.
Note that $\GD$ is as follows 
\begin{alignat*}{2}
\langle \diag(a,b,c)\colon &    a,b,c\in \FF_{q^2}\mbox{ with } aa^\sigma=bb^\sigma=cc^\sigma=1 \rangle && \mbox{ if } \liediag=A_2,\\
\langle \diag(c^\sigma b, ab ,c,a^\sigma)\colon &   a,b,c\in \FF_{q^2}\mbox{ with } aa^\sigma=bb^\sigma=cc^\sigma=1 \rangle && \mbox{ if } \liediag=C_2.
\end{alignat*}
\ele

\bre
\begin{enumerate}
\item In case $\liediag=C_2$, $\amgrpG\cong\Sp_4(q)$ is realized as $\Sp_4(q^2)\cap \SU_4(q)$
 with respect to a basis that is 
 hyperbolic for the symplectic form and orthonormal for the unitary form, and $\Aut(\FF_{q^2})$ acts entry-wise on these matrices. 
Moreover, $\tau$ acts as transpose-inverse on these matrices.
\item In all cases $\trin$ coincides with $\sigma$. 
\end{enumerate}
\ere

\bpf
The $A_1\times A_1$ case is self evident.
Now consider the case $\liediag=A_2$.
As in the proof of Lemma~\ref{lem:structure of Phan A_i groups}, $\Aut(\FF_{q^2})\le N_{\Aut(\amgrpG_{i,j})}(\amgrpG_i)\cap N_{\Aut(\amgrpG_{i,j})}(\amgrpG_j)$, $A^\trin={}^tA^{-1}=A^\sigma$ and $\Aut(\amgrpG_{i,j})\cong\PGU_3(q)\rtimes\Aut(\FF_{q^2})$, so it suffices to consider linear automorphisms.
As before this is an uncomplicated calculation.

Now consider the case $\liediag=C_2$. Writing $\GamL(V)\cong \GL_4(q^2)\rtimes\Aut(\FF_{q^2})$ with respect to the basis $\cE=\{e_1,e_2,e_3=f_1,e_4=f_2\}$, which is hyperbolic for the symplectic form of $\Sp_4(q^2)$ and orthonormal for the unitary form of $\SU_4(q)$, we have $\amgrpG_{i,j}=\Sp_4(q^2)\cap \SU_4(q)$.

There is an isomorphism $\Phi\colon \amgrpG_{i,j}\to \Sp_4(q)$ as in~\cite{GraHofShp2003}.
Abstractly, we have  $\Aut(\Sp_4(q))= \GSp_4(q)\rtimes\Aut(\FF_q)$ (with respect to a suitable basis $\sfE$ for $V$).
Since the embedding of $\Sp_4(q)$ into $\Sp_4(q^2)$ is non-standard, we are reconstructing the automorphism group here.

We first note that changing bases just replaces $\Aut(\FF_{q^2})$ with a different complement to the linear automorphism group.
As for linear automorphisms we claim that 
\begin{align*}
\GSp_4(q^2)\cap \GU_4(q)=\GSp_4(q)
\end{align*}
 (viewing the latter as a matrix group w.r.t.~$\sfE$).
Clearly, up to a center, we have  $\amgrpG_{i,j}\le \GSp_4(q^2)\cap \GU_4(q)\le \GSp_4(q)$ and we note that 
 $\GSp_4(q)/\Sp_4(q)\cong (\FF_q^*)^2/ (\FF_q^*)$.
Thus for $q$ even, the claim follows.
For $q$ odd, let $\FF_{q^2}^*=\langle \zeta\rangle$ and define  $\beta=\diag(\zeta^{q-1},\zeta^{q-1},1,1)$.
Then $\beta\in \GSp_4(q^2)\cap \GU_4(q)$ acts  
 on $\amgrpG_{i,j}$ as $\diag(\zeta^q,\zeta^q,\zeta,\zeta)$, which scales the symplectic form of $\Sp_4(q^2)$ by $\zeta^{q+1}$.
By~\cite{GraHofShp2003} the form of $\Sp_4(q)$ is proportional and since $\zeta^{q+1}$ is a non-square in $\FF_q$,  $\beta$ is a linear outer automorphism of $\Sp_4(q)$.  Thus, $\GSp_4(q)=\langle \Sp_4(q),\beta\rangle$ and the claim follows.

We now determine $\amgrpA_{i,j}$.
First we note that $\beta$, as well as the group $\Aut(\FF_{q^2})$ with respect to the basis $\cE$, clearly normalize $\amgrpG_i$ and $\amgrpG_j$ hence by Lemma~\ref{lem:N G1 G2=N D1 D2}, $\Aut(\FF_{q^2})\le \amgrpA_{i,j}$. So it suffices to determine inner automorphisms of $\Sp_4(q)$ normalizing $\amgrpD_i^j$ and $\amgrpD_j^i$.

Any inner automorphism in $\Sp_4(q)$ is induced by an inner automorphism of $\Sp_4(q^2)$.
So now the claim reduces to a matrix calculation in the group $\Sp_4(q^2)$.
\epf

\medskip

Next we describe the restriction maps $\ama_{i,j}$ for Phan amalgams made up of a single standard pair with trivial inclusion maps.

\ble\label{lem: Phan coefficient system connecting maps}
Let $\amA$ be the coefficient system of the standard Phan  amalgam $\famG$ over $\FF_q$ with diagram $\liediag$ and system of tori $\cD$. 
Fix $i,j\in I$ and let $(\amgrpG_{i,j},\bamgrpG_i,\bamgrpG_j)$ be a Phan standard pair in $\famG$ with diagram $\liediag_{i,j}$. Denote $\ama=(\ama_{j,i},\ama_{i,j})\colon \amgrpA_{i,j}\to \amgrpA_i\times\amgrpA_j$.
Then, we have the following:
\begin{enumerate}
\item If $\liediag_{i,j}=A_1\times A_1$, then $\ama$ is an isomorphism inducing $\amgrpT_{i,j}\cong\amgrpT_i\times\amgrpT_j$ and $\amgrpC_{i,j}\cong\amgrpC_i\times\amgrpC_j$.
\item If $\liediag_{i,j}=A_2$ or $C_2$, then $\ama$ induces an isomorphism $\amgrpT_{i,j}\to \amgrpT_i\times\amgrpT_j$.
\item If $\liediag(q)=A_2(q)$ or $\liediag(q)=C_2(q)$, then $\ama\colon\amgrpC_{i,j}\to \amgrpC_i\times\amgrpC_j$ is given by 
  $\alpha\mapsto (\alpha,\alpha)$ (for $\alpha\in \Aut(\FF_{q^2})$) which is a diagonal  embedding.
\end{enumerate}
\ele
\bpf
1.~ This is immediate from Lemma~\ref{lem:structure of Phan A_ij groups}.

For the remaining cases, recall that for any $\varphi\in \amgrpA_{i,j}$, we have 
 $\ama_{i,j}\colon \varphi\mapsto \famg_{j,i}^{-1}\after \rho_{i,j}(\varphi)\after\famg_{j,i}$, 
where $\rho_{i,j}(\varphi)$ is the restriction of $\varphi$ to $\bamgrpG_j\le \amgrpG_{i,j}$ (Definition~\ref{dfn:phan coefficient system}) and $\famg_{i,j}$ is the standard identification map of Definition~\ref{dfn:standard Phan identification map}.
Note that the standard identification map transforms the automorphism $\rho_{j,i}(\varphi)$ of $\bamgrpG_i$ essentially to the ``same" automorphism $\varphi$ of $\amgrpG_i$.

First let $\liediag(q)=A_2(q)$.
The map $\ama$ is well-defined. On $\GD$, it is induced by the homomorphism
 \begin{align*}
\diag(ac,c,ec) \mapsto (\diag(1,e), \diag(a,1)), 
 \end{align*}
where $a,c,e\in \FF_{q^2}$ are such that $aa^\sigma=cc^\sigma=ee^\sigma=1$.
Note that the kernel is $Z(\GD)$ so that $\ama$ is injective. The map is obviously surjective, so we are done.
Thus if we factor $\ama$ by $\amgrpT_{i,j}$ and $\amgrpT_i\times\amgrpT_j$, we get
\begin{align}
\amgrpC_{i,j}  \into &  \amgrpC_i\times\amgrpC_j \label{eqn:rho mod PGD}
 \end{align}
which is a diagonal embedding given by $\alpha^r\mapsto (\alpha^r,\alpha^r)$, where 
  $r\in \NN$ and $\alpha\colon x\mapsto x^p$ for $x\in \FF_{q^2}$.

Next let $\liediag(q)=C_2(q)$.
We can rewrite the elements of $\GD$ as 
 a diagonal matrix $\diag(xyz, xz, z y^{-1}, z)$, by taking $z=a^{\sigma}$, $y=(ac)^{-1}$, $x=a^2b$.
On $\GD$ the map $\ama$ is induced by the homomorphism 
 \begin{align*}
\diag(xyz,xz,y^{-1}z,z) \mapsto (\diag(y,1),\diag(x,1))
 \end{align*}
with kernel $\{\diag(z,z,z,z)\colon z\in \FF_{q^2} \mbox{ with }zz^\sigma=1\}=Z(\GU_4(q))$. 
Clearly $\ama\colon\amgrpT_{i,j}\to \amgrpT_i\times\amgrpT_j$ is an isomorphism.
Taking the quotient over these groups, $\ama$ induces a diagonal embedding as in~\eqref{eqn:rho mod PGD}, where we now interpret it in the $C_2(q)$ setting.
\epf

\medskip
\subsubsection{A standard form for Phan amalgams}\label{subsec:Phan trivial support}

Suppose that $\famG=\{\amgrpG_i,\amgrpG_{i,j},\famg_{i,j}\mid i,j\in I\}$ is a Phan amalgam over $\FF_q$ with $3$-spherical diagram $\liediag$.
Without loss of generality we will assume that all inclusion maps $\famg_{i,j}$ are the standard identification maps of Definition~\ref{dfn:standard Phan identification map}.
By Lemma~\ref{lem:non-collapsing Phan has property D}, $\famG$ has Property (D) and possesses a system 
 $\cD=\{\amgrpD_i\colon i\in I\}$ of tori, which, as noted in Definition~\ref{dfn:type of P amalgam} via the standard embeddings $\famg_{i,j}$ can be identified with those given in Lemma~\ref{lem:tori in phan standard pairs}.

We wish to classify all  Phan amalgams  $\amG=\{\amgrpG_i,\amgrpG_{i,j},\amg_{i,j}\mid i,j\in I\}$ over $\FF_q$ with the same diagram as $\famG$.
As noted in Definition~\ref{dfn:type of P amalgam} we may assume that all such amalgams share $\cD$.
Let $\amA=\{\amgrpA_i,\amgrpA_{i,j},\ama_{i,j}\mid i,j\in I\}$ be the coefficient system of $\famG$ associated to $\cD$.
By Proposition~\ref{prop:P coefficient systems}, we may restrict to those amalgams whose connecting maps 
 are of the form $\amg_{i,j}=\famg_{i,j}\after \delta_{i,j}$ for $\delta_{i,j}\in \amgrpA_i$ for all $i\in I$.

\bde\label{dfn:Phan trivial support}
The {\dfn trivial support} of $\amG$ (with respect to $\famG$) is the set $\{(i,j)\in I\times I\mid \amg_{i,j}=\famg_{i,j}\}$ (that is, $\delta_{i,j}=\id_{\amgrpG_i}$ in the notation of Proposition~\ref{prop:P coefficient systems}).
The word ``trivial'' derives from the assumption that the $\famg_{i,j}$'s are the standard identification maps of Definition~\ref{dfn:standard Phan identification map}.
\ede

Fix some spanning tree $\Sigma\sbe \Gamma$ and suppose that $\edg-\edg\Sigma=\{\{i_s,j_s\}\colon s=1,2,\ldots,r\}$
 so that $H_1(\Gamma,\ZZ)\cong\ZZ^r$. 
We now have 
\bpr\label{prop:Phan trivial support on spanning tree}
There is a  Phan amalgam $\amG(\Sigma)$ with the same diagram as $\famG$ and the same $\cD$, which is isomorphic to $\amG$ and has the following properties:
\begin{enumerate}
\item $\amG$ has trivial support $S=\{(i,j)\in I\times I\mid \{i,j\}\in \edg \Sigma\}\cup \{(i_s,j_s)\colon s=1,2\ldots,r\}$.
\item for each  $s=1,2,\ldots,r$, we have $\amg_{j_s,i_s}=\famg_{j_s,i_s}\after\gamma_{j_s,i_s}$, where $\gamma_{j_s,i_s}\in \amgrpC_{j_s}$.
\end{enumerate}
\epr

\ble\label{lem:Phan stripping tori from connecting maps}
There is a  Phan  amalgam $\amG^+$ over $\FF_q$ with the same diagram as $\famG$ and the same $\cD$, which is isomorphic to $\amG$ and has the following properties:
For any $u,v\in I$, if $\amg_{u,v}=\famg_{u,v}\after \gamma_{u,v}\after d_{u,v}$,  for some $\gamma_{u,v}\in \amgrpC_{u}$ and $d_{u,v}\in \amgrpT_u$, then $\amg_{u,v}^+=\famg_{u,v}\after\gamma_{u,v}$.
\ele
\bpf
The proof follows the same steps as that of Lemma~\ref{lem:stripping tori from connecting maps} using Part 2 of Lemma~\ref{lem: Phan coefficient system connecting maps} instead of Lemma~\ref{lem: coefficient system connecting maps} Part 2. 
\epf

\medskip
By Lemma~\ref{lem:Phan stripping tori from connecting maps} in order to prove Proposition~\ref{prop:Phan trivial support on spanning tree} we may now assume that $\amg_{u,v}=\famg_{u,v}\after\gamma_{u,v}$ for some $\gamma_{u,v}\in \amgrpC_u$ for all $u,v\in I$.

We now prove a Corollary for Phan amalgams analogous to, but stronger than Corollary~\ref{cor:completing the hexagon}. To this end consider the situation of Figure~\ref{fig:completing the hexagon} interpreted in the Phan setting. 
\bco\label{cor:Phan completing the hexagon}
With the notation introduced in Figure~\ref{fig:completing the hexagon}, fix the maps $\gamma_{i,j}$, $\gamma^+_{i,j}$, $\phi_i \in \amgrpC_i$ as well as 
 $\gamma_{j,i}\in \amgrpC_j$. 
Then for any one of $\gamma^+_{j,i}, \phi_j\in \amgrpC_j$, there exists a unique choice $\gamma\in \amgrpC_i$ for the remaining map in $\amgrpC_j$  so that there exists $\phi_{i,j}$ making the diagram in Figure~\ref{fig:completing the hexagon} commute.
\eco
\bpf
This  follows immediately from the fact that the maps $\ama_{j,i}\colon \amgrpC_{i,j}\to \amgrpC_i$ and $\ama_{i,j}\colon \amgrpC_{i,j}\to \amgrpC_j$ in part 3.~of Lemma~\ref{lem: Phan coefficient system connecting maps} are isomorphisms.
\epf

\medskip
\bpf (of Proposition~\ref{prop:Phan trivial support on spanning tree})
The proof follows the same steps as that of Proposition~\ref{prop:trivial support on spanning tree}, replacing 
Lemma~\ref{lem:stripping tori from connecting maps}~and Corollary~\ref{cor:completing the hexagon} by 
Lemma~\ref{lem:Phan stripping tori from connecting maps}~and Corollary~\ref{cor:Phan completing the hexagon}.
\epf

\medskip
\subsubsection{Classification of  Phan amalgams with $3$-spherical diagram}\label{subsubsec:classification of 3-spherical Phan amalgams}
In the case where $\amG$ is a Phan  amalgam over $\FF_q$ whose diagram is a $3$-spherical tree,  Proposition~\ref{prop:Phan trivial support on spanning tree} says that $\amG\cong\famG$. 

\bth\label{thm:Phan 3 spherical tree}
Suppose that $\amG$ is a  Curtis-Tits amalgam with a diagram that is a $3$-spherical tree.
Then, $\amG$ is unique up to isomorphism. In particular any Phan amalgam with spherical diagram is unique.
\eth

\bde\label{dfn:P the map kappa}
Fix a connected $3$-spherical diagram $\liediag$ and a prime power $q$. Let  $\Sigma$  be a spanning tree and let the set of edges $\edg\liediag-\edg\Sigma=\{\{i_s,j_s\}\colon s=1,2,\ldots,r\}$ together with the integers $\{e_s\colon s=1,2,\ldots,r\}$ satisfy the conclusions of Lemma~\ref{lem:spanning tree with minimal complement}.
Note that since in the Phan case we do not have subdiagrams of type $\twA_3(q)$, we have $e_s=1$ for all $s=\{1,2,\ldots,r\}$.

Let $\sPh(\liediag,q)$ be the collection of isomorphism classes of  Phan  amalgams of type $\liediag(q)$ and let $\famG=\{\amgrpG_i,\amgrpG_{i,j},\famg_{i,j}\mid i,j\in I\}$ be a  Phan amalgam over $\FF_q$ with diagram $\liediag$.

Consider the following map:
\begin{align*}
\kappa\colon\prod_{s=1}^r \Aut(\FF_{q^2})\to \sPh(\liediag).
\end{align*}
where $\kappa((\alpha_s)_{s=1}^r)$ is the isomorphism class of the amalgam $\amG^+=\amG((\alpha_s)_{s=1}^r)$ given by setting $\amg^+_{j_s,i_s}=\amg_{j_s,i_s}\after \alpha_s$ for all $s=1,2,\ldots,r$.
\ede

As for Curtis-Tits amalgams, one shows the following.
\bco\label{cor:P kappa is onto}
The map $\kappa$ is onto.
\eco

\ble\label{lem:classification of  P amalgams with loop diagram}
Suppose $\liediag(q)$ is a 
$3$-spherical diagram $\liediag$ that is a simple loop. 
Then, $\kappa$ is injective.
\ele

\bpf
The proof is identical to that of Lemma~\ref{lem:classification of  CT amalgams with loop diagram} replacing Proposition~\ref{prop:trivial support on spanning tree} by Proposition~\ref{prop:Phan trivial support on spanning tree}  and Lemma~\ref{lem: coefficient system connecting maps} by Lemma~\ref{lem: Phan coefficient system connecting maps}, and noting that in the Phan case, we can consider the group $\amgrpC_i$ and $\amgrpC_{i,j}$ themselves rather than some suitably chosen quotient.
\epf

\bth\label{thm:P classification of 3-spherical amalgams}
Let  $\liediag$ be a connected $3$-spherical diagram with spanning tree $\Sigma$ and set of edges
 $\edg\liediag-\edg\Sigma=\{\{i_s,j_s\}\colon s=1,2,\ldots,r\}$.
Then $\kappa$ is a bijection between the elements of $\prod_{s=1}^r \Aut(\FF_{q^2})$ and the isomorphism classes of  Curtis-Tits amalgams with diagram $\liediag$ over $\FF_q$.
\eth
\bpf
This follows from Lemma~\ref{lem:classification of  P amalgams with loop diagram}
just as Theorem~\ref{thm:CT classification of 3-spherical amalgams}
follows from Lemma~\ref{lem:classification of  CT amalgams with loop diagram}.
\epf

\bibliographystyle{alpha}

\begin{thebibliography}{10}

%\bibitem{AbrMuh97}
%P.~Abramenko and B.~M{\"u}hlherr.
%\newblock Pr\'esentations de certaines {$BN$}-paires jumel\'ees comme sommes
%  amalgam\'ees.
%\newblock {\em C. R. Acad. Sci. Paris S\'er. I Math.}, 325(7):701--706, 1997.
%
%\bibitem{BeGrHoSh2003}
%C.~D. Bennett, R.~Gramlich, C.~Hoffman, and S.~Shpectorov.
%\newblock Curtis-{P}han-{T}its theory.
%\newblock In {\em Groups, combinatorics \& geometry ({D}urham, 2001)}, pages
%  13--29. World Sci. Publ., River Edge, NJ, 2003.
%
%\bibitem{BenGraHof2003}
%C.~D. Bennett, R.~Gramlich, C.~Hoffman, and S.~Shpectorov.
%\newblock Curtis-{P}han-{T}its theory.
%\newblock In {\em Groups, combinatorics \& geometry ({D}urham, 2001)}, pages
%  13--29. World Sci. Publ., River Edge, NJ, 2003.

\bibitem{BeSh2004}
C.~D. Bennett and S.~Shpectorov.
\newblock A new proof of a theorem of {P}han.
\newblock {\em J. Group Theory}, 7(3):287--310, 2004.

%\bibitem{BlHo2008}
%R.~J. Blok and C.~Hoffman.
%\newblock A quasi {C}urtis-{T}its-{P}han theorem for the symplectic group.
%\newblock {\em J. Algebra}, 319(11):4662--4691, 2008.
%
%\bibitem{BlHo2009}
%R.~J. Blok and C.~Hoffman.
%\newblock A {C}urtis-{T}its-{P}han theorem for the twin-building of type
%  {$\widetilde A\sb {n-1}$}.
%\newblock {\em J. Algebra}, 321(4):1196--1124, 2009.
%
%\bibitem{BloHof2011}
%R.~J. Blok and C.~G.~Hoffman.
%\newblock Bass-{S}erre theory and counting rank two amalgams.
%\newblock {\em J. Group Theory}, 14(3):389--400, 2011.
%
\bibitem{BloHof2013}
R.~J. Blok and C.~G. Hoffman.
\newblock 1-cohomology of simplicial amalgams of groups.
\newblock {\em J. Algebraic Combin.}, 37(2):381--400, 2013.
%
\bibitem{BloHof2014a}
R.~J. Blok and C.~G. Hoffman.
\newblock Curtis--{T}its groups generalizing {K}ac--{M}oody groups of type
  {$\tilde{A}_{n-1}$}.
\newblock {\em J. Algebra}, 399:978--1012, 2014.

\bibitem{BloHof2014b}
R.~J. Blok and C.~Hoffman.
\newblock A classfication of {C}urtis-{T}its amalgams.
\newblock In N.~Sastry, editor, {\em Groups of Exceptional Type, Coxeter Groups
  and Related Geometries}, volume 149 of {\em Springer Proceedings in
  Mathematics \& Statistics}, pages 1--26. Springer, January 2014.
%
\bibitem{BloHof2016}
R.~J. Blok and C.~G. Hoffman.
\newblock Curtis-Tits groups of simply-laced type.
\newblock To appear in J.~Comb.~Th.~Ser. A.
%
%
%\bibitem{Bray:2013aa}
%J.~N. Bray, D.~F. Holt, and C.~M. Roney-Dougal.
%\newblock {\em The maximal subgroups of the low-dimensional finite classical
%  groups}, volume 407 of {\em London Mathematical Society Lecture Note Series}.
%\newblock Cambridge University Press, Cambridge, 2013.
%\newblock With a foreword by Martin Liebeck.
%

%\bibitem{BloHof2014a}
%R.~J. Blok and C.~G. Hoffman.
%\newblock Curtis--{T}its groups generalizing {K}ac--{M}oody groups of type
%  {$\tilde{A}_{n-1}$}.
%\newblock {\em J. Algebra}, 399:978--1012, 2014.

%\bibitem{BloHof2014b}
%R.~J.~Blok and C.~G.~Hoffman.
%\newblock A classification of {C}urtis-{T}its amalgams.
%\newblock To appear in {\em Groups of Exceptional Type, Coxeter Groups and Related Geometries}, 
%Springer Proceedings in Mathematics \& Statistics, Vol. 149, Sastry, N.S. Narasimha (Ed.).
%
\bibitem{Cap2007}
P.~E.~Caprace.
\newblock On 2-spherical {K}ac-{M}oody groups and their central extensions.
\newblock {\em Forum Math.}, 19(5):763--781, 2007.

%\bibitem{Ca1972}
%R.~W.~Carter.
%\newblock {\em Simple groups of Lie type}, volume~28 of {\em Pure and Applied
%  Math.}
%\newblock Wiley, London, 1972.

\bibitem{Cur1965a}
C.~W.~Curtis.
\newblock Central extensions of groups of {L}ie type.
\newblock {\em J. Reine Angew. Math.}, 220:174--185, 1965.

%\bibitem{DevMuh2007}
%A.~Devillers and B.~M{\"u}hlherr.
%\newblock On the simple connectedness of certain subsets of buildings.
%\newblock {\em Forum Math.}, 19(6):955--970, 2007.

\bibitem{Dun2005}
J.~Dunlap.
\newblock {\em Uniqueness of Curtis-Phan-Tits amalgams}.
\newblock PhD thesis, Bowling Green State University, 2005.
%
%
%\bibitem{GorLyoSol1998}
%D.~Gorenstein, R.~Lyons, and R.~Solomon.
%\newblock {\em The classification of the finite simple groups. {N}umber 3.
%  {P}art {I}. {C}hapter {A}}, volume~40 of {\em Mathematical Surveys and
%  Monographs}.
%\newblock American Mathematical Society, Providence, RI, 1998.
%\newblock Almost simple $K$-groups.
%
%
%\bibitem{GraHorMuh2011}
%R.~Gramlich, M.~Horn, and B.~M{\"u}hlherr.
%\newblock Abstract involutions of algebraic groups and of {K}ac-{M}oody groups.
%\newblock {\em J. Group Theory}, 14(2):213--249, 2011.
\bibitem{Gor1983}
D.~Gorenstein.
\newblock {\em The classification of finite simple groups. {V}ol. 1}.
\newblock The University Series in Mathematics. Plenum Press, New York, 1983.
\newblock Groups of noncharacteristic $2$ type.

\bibitem{GorLyoSol1996}
D.~Gorenstein, R.~Lyons, and R.~Solomon.
\newblock {\em The classification of the finite simple groups. {N}umber 2.
  {P}art {I}. {C}hapter {G}}, volume~40 of {\em Mathematical Surveys and
  Monographs}.
\newblock American Mathematical Society, Providence, RI, 1996.
\newblock General group theory.

\bibitem{GorLyoSol1998}
D.~Gorenstein, R.~Lyons, and R.~Solomon.
\newblock {\em The classification of the finite simple groups. {N}umber 3.
  {P}art {I}. {C}hapter {A}}, volume~40 of {\em Mathematical Surveys and
  Monographs}.
\newblock American Mathematical Society, Providence, RI, 1998.
\newblock Almost simple $K$-groups.

\bibitem{GorLyoSol1999}
D.~Gorenstein, R.~Lyons, and R.~Solomon.
\newblock {\em The classification of the finite simple groups. {N}umber 4.
  {P}art {II}. {C}hapters 1--4}, volume~40 of {\em Mathematical Surveys and
  Monographs}.
\newblock American Mathematical Society, Providence, RI, 1999.
\newblock Uniqueness theorems, With errata: {{\i}t The classification of the
  finite simple groups. Number 3. Part I. Chapter A} [Amer. Math. Soc.,
  Providence, RI, 1998; MR1490581 (98j:20011)].

\bibitem{GorLyoSol2002}
D.~Gorenstein, R.~Lyons, and R.~Solomon.
\newblock {\em The classification of the finite simple groups. {N}umber 5.
  {P}art {III}. {C}hapters 1--6}, volume~40 of {\em Mathematical Surveys and
  Monographs}.
\newblock American Mathematical Society, Providence, RI, 2002.
\newblock The generic case, stages 1--3a.

\bibitem{GorLyoSol2005}
D.~Gorenstein, R.~Lyons, and R.~Solomon.
\newblock {\em The classification of the finite simple groups. {N}umber 6.
  {P}art {IV}}, volume~40 of {\em Mathematical Surveys and Monographs}.
\newblock American Mathematical Society, Providence, RI, 2005.
\newblock The special odd case.

\bibitem{Gr2004}
R.~Gramlich.
\newblock Weak {P}han systems of type {$C\sb n$}.
\newblock {\em J. Algebra}, 280(1):1--19, 2004.



%\bibitem{Mu1993}
%B.~M{\"u}hlherr.
%\newblock Coxeter groups in {C}oxeter groups.
%\newblock In {\em Finite geometry and combinatorics (Deinze, 1992)}, volume 191
%  of {\em London Math. Soc. Lecture Note Ser.}, pages 277--287. Cambridge Univ.
%  Press, Cambridge, 1993.


%\bibitem{Re2002a}
%B.~R{\'e}my.
%\newblock Groupes de {K}ac-{M}oody d\'eploy\'es et presque d\'eploy\'es.
%\newblock {\em Ast\'erisque}, (277):viii+348, 2002.
%
%
%\bibitem{Re2004}
%B.~R{\'e}my.
%\newblock Kac-{M}oody groups: split and relative theories. {L}attices.
%\newblock In {\em Groups: topological, combinatorial and arithmetic aspects},
%  volume 311 of {\em London Math. Soc. Lecture Note Ser.}, pages 487--541.
%  Cambridge Univ. Press, Cambridge, 2004.

%\bibitem{Spr1985}
%T.~A.~Springer.
%\newblock Some results on algebraic groups with involutions.
%\newblock In {\em Algebraic groups and related topics ({K}yoto/{N}agoya,
%  1983)}, volume~6 of {\em Adv. Stud. Pure Math.}, pages 525--543.
%  North-Holland, Amsterdam, 1985.

%\bibitem{St1967}
%R.~Steinberg.
%\newblock {\em Lectures on Chevalley groups}.
%\newblock Yale Lecture Notes. Yale University, 1967.

\bibitem{Tim1998}
F.~G.~Timmesfeld.
\newblock Presentations for certain {C}hevalley groups.
\newblock {\em Geom. Dedicata}, 73(1):85--117, 1998.

\bibitem{Tim03}
F.~G.~Timmesfeld.
\newblock On the {S}teinberg-presentation for {L}ie-type groups.
\newblock {\em Forum Math.}, 15(5):645--663, 2003.

\bibitem{Tim04}
F.~G.~Timmesfeld.
\newblock The {C}urtis-{T}its-presentation.
\newblock {\em Adv. Math.}, 189(1):38--67, 2004.

\bibitem{Tim06}
F.~G.~Timmesfeld.
\newblock Steinberg-type presentation for {L}ie-type groups.
\newblock {\em J. Algebra}, 300(2):806--819, 2006.

%\bibitem{Ti1974}
%J.~Tits.
%\newblock {\em Buildings of spherical type and finite {BN}-pairs}.
%\newblock Springer-Verlag, Berlin, 1974.
%\newblock Lecture Notes in Mathematics, Vol. 386.
%
%\bibitem{Ti1986b}
%J.~Tits.
%\newblock Ensembles Ordonn\'es, immeubles et sommes amalgam\'ees.
%\newblock {\em Bull. Soc. Math. Belg. S\'er A} 38:367-387, 1986.
%
%\bibitem{Ti1987}
%J.~Tits.
%\newblock Uniqueness and presentation of {K}ac-{M}oody groups over fields.
%\newblock {\em J. Algebra}, 105(2):542--573, 1987.
%
%
%\bibitem{Ti1992}
%J.~Tits.
%\newblock Twin buildings and groups of {K}ac-{M}oody type.
%\newblock In {\em Groups, combinatorics \& geometry (Durham, 1990)}, volume 165
%  of {\em London Math. Soc. Lecture Note Ser.}, pages 249--286. Cambridge Univ.
%  Press, Cambridge, 1992.



%\bibitem{BloHof2013}
%R.~J. Blok and C.~G. Hoffman.
%\newblock 1-cohomology of simplicial amalgams of groups.
%\newblock {\em J. Algebraic Combin.}, 37(2):381--400, 2013.


%\bibitem{FiSo1979}
%L.~Finkelstein and R.~Solomon.
%\newblock A presentation of the symplectic and orthogonal groups.
%\newblock {\em J. Algebra}, 60(2):423--438, 1979.

\bibitem{Gra2004}
R.~Gramlich.
\newblock Weak {P}han systems of type {$C_n$}.
\newblock {\em J. Algebra}, 280(1):1--19, 2004.
%
%\bibitem{Gr2004}
%R.~Gramlich.
%\newblock Weak {P}han systems of type {$C\sb n$}.
%\newblock {\em J. Algebra}, 280(1):1--19, 2004.

%\bibitem{Gra2009}
%R.~Gramlich.
%\newblock Developments in finite {P}han theory.
%\newblock {\em Innov. Incidence Geom.}, 9:123--175, 2009.

\bibitem{GraHofShp2003}
R.~Gramlich, C.~Hoffman, and S.~Shpectorov.
\newblock A {P}han-type theorem for {${\rm Sp}(2n,q)$}.
\newblock {\em J. Algebra}, 264(2):358--384, 2003.

%\bibitem{GrHoSh2003}
%R.~Gramlich, C.~Hoffman, and S.~Shpectorov.
%\newblock A {P}han-type theorem for {${\rm Sp}(2n,q)$}.
%\newblock {\em J. Algebra}, 264(2):358--384, 2003.
%
%\bibitem{GrHoSh2003a}
%R.~Gramlich, C.~Hoffman, and S.~Shpectorov.
%\newblock A {P}han-type theorem for {${\rm Sp}(2n,q)$}.
%\newblock {\em J. Algebra}, 264(2):358--384, 2003.

\bibitem{GraHorNic2006}
R.~Gramlich, M.~Horn, and W.~Nickel.
\newblock The complete {P}han-type theorem for {${\rm Sp}(2n,q)$}.
\newblock {\em J. Group Theory}, 9(5):603--626, 2006.

%\bibitem{GraHorNic2006a}
%R.~Gramlich, M.~Horn, and W.~Nickel.
%\newblock The complete {P}han-type theorem for {${\rm Sp}(2n,q)$}.
%\newblock {\em J. Group Theory}, 9(5):603--626, 2006.

\bibitem{GrHoNi2006}
R.~Gramlich, M.~Horn, and W.~Nickel.
\newblock The complete {P}han-type theorem for {${\rm Sp}(2n,q)$}.
\newblock {\em J. Group Theory}, 9(5):603--626, 2006.

%\bibitem{Hoffman:2013aa}
%C.~Hoffman and A.~Roberts.
%\newblock On a quasi-{P}han theorem for orthogonal groups.
%\newblock {\em Comm. Algebra}, 41(5):1589--1600, 2013.
%
%\bibitem{Hor2008}
%M.~Horn.
%\newblock On the {P}han system of the {S}chur cover of {${\rm SU}(4,3^2)$}.
%\newblock {\em Des. Codes Cryptogr.}, 47(1-3):243--247, 2008.
%
%
%\bibitem{KleLie1990a}
%P.~Kleidman and M.~Liebeck.
%\newblock {\em The subgroup structure of the finite classical groups}, volume
%  129 of {\em London Mathematical Society Lecture Note Series}.
%\newblock Cambridge University Press, Cambridge, 1990.
%
\bibitem{Mu1999}
B.~M{\"u}hlherr.
\newblock Locally split and locally finite twin buildings of {$2$}-spherical
  type.
\newblock {\em J. Reine Angew. Math.}, 511:119--143, 1999.

%\bibitem{MuRo1995}
%B.~M{\"u}hlherr and M.~Ronan.
%\newblock Local to global structure in twin buildings.
%\newblock {\em Invent. Math.}, 122(1):71--81, 1995.
%
%
%
\bibitem{Pha1971}
K.-W. Phan.
\newblock A characterization of the unitary groups {${\rm
  PSU}(4,\,q^{2}),\,q$}\ odd.
\newblock {\em J. Algebra}, 17:132--148, 1971.

\bibitem{Pha1977}
K.~W. Phan.
\newblock On groups generated by three-dimensional special unitary groups. {I}.
\newblock {\em J. Austral. Math. Soc. Ser. A}, 23(1):67--77, 1977.

\bibitem{Pha1977a}
K.-W. Phan.
\newblock On groups generated by three-dimensional special unitary groups.
  {II}.
\newblock {\em J. Austral. Math. Soc. Ser. A}, 23(2):129--146, 1977.
%
%\bibitem{Ro1989a}
%M.~Ronan.
%\newblock {\em Lectures on buildings}, volume~7 of {\em Perspectives in
%  Mathematics}.
%\newblock Academic Press Inc., Boston, MA, 1989.
%
\bibitem{SchVan28}
O.~Schreier and B.~Van~der Waerden.
\newblock Die automorphismen der projektiven gruppen.
\newblock {\em Abhandlungen aus dem Mathematischen Seminar der Universit\"{a}t
  Hamburg}, 6(1):303--322, December 1928.
%
\bibitem{Ti1992}
J.~Tits.
\newblock Twin buildings and groups of {K}ac-{M}oody type.
\newblock In {\em Groups, combinatorics \& geometry (Durham, 1990)}, volume 165
  of {\em London Math. Soc. Lecture Note Ser.}, pages 249--286. Cambridge Univ.
  Press, Cambridge, 1992.
%
\bibitem{Wil2009}
R.~A.~Wilson.
\newblock {\em The finite simple groups}, volume 251 of {\em Graduate Texts in
  Mathematics}.
\newblock Springer-Verlag London Ltd., London, 2009.
%

\end{thebibliography}

\end{document}